\documentclass[12pt,letterpaper]{article}
\usepackage[latin1]{inputenc}
\usepackage{amsfonts}
\usepackage{amsmath}
\usepackage{amstext}
\usepackage{amsbsy}
\usepackage{amssymb}
\usepackage{graphics}
\usepackage{graphicx}
\usepackage{yfonts}
\usepackage{pstricks-add}
\usepackage{pst-math}
\usepackage{url}

\usepackage{authblk}

\newtheorem{thm}{Theorem}[section]
\newcommand{\bthm}{\begin{thm}\noindent{\bf }~}
\newcommand{\ethm}{\end{thm}}
\newtheorem{propo}{Proposition}[section]
\newcommand{\bprop}{\begin{propo}\noindent{\bf }~}
\newcommand{\eprop}{\end{propo}}
\newtheorem{lema}{Lemma}[section]
\newcommand{\blem}{\begin{lema}\noindent{\bf }~}
\newcommand{\elem}{\end{lema}}
\newtheorem{defe}{Definition}[section]
\newcommand{\bdefe}{\begin{defe}\noindent{\bf }~}
\newcommand{\edefe}{\end{defe}}
\newtheorem{ejem}{Example}[section]
\newcommand{\bejem}{\begin{ejem}\noindent{\bf }~}
\newcommand{\eejem}{\end{ejem}}
\newtheorem{coro}{Corollary}[section]
\newcommand{\bcor}{\begin{coro}\noindent{\bf }~}
\newcommand{\ecor}{\end{coro}}
\newtheorem{rema}{Remark}[section]
\newcommand{\brem}{\begin{rema}\noindent{\bf }~\rm}
\newcommand{\erem}{\end{rema}}

\newcommand{\edemo}{\hfill{$\sqcap \kern-18pt \sqcup$}}

\title{\large\bf Geodesibility of algebrizable three-dimensional vector fields}

\author[D. Kopta et al.]
       {{\small  M. E. Fr\'ias-Armenta$^1$,  E. L\'opez-Gonz\'alez$^2$}
       \\
       {\footnotesize{$^1$Universidad de Sonora, Departamento de Matem\'aticas, Blvrd. Rosales
y Luis Encinas S/N,
Col. Centro, Hermosillo Sonora, C.P. 83000, M\'exico.\\
\ttfamily eduardo@gauss.mat.uson.mx\\
       $^2$Universidad Aut\'onoma de Ciudad Ju\'arez,
Unidad Multidisciplinaria de la UACJ en Cuauht\'emoc, Carretera
Cuauht\'emoc-An\'ahuac, Col. Ejido
An\'ahuac Km. 3.5 S/N, Mpio. de Cuauht\'emoc, Chih., C.P. 31600, M\'exico.\\
       \ttfamily elgonzal@uacj.mx\\
       }}}

\begin{document}

\maketitle

\begin{small}\noindent \textbf{Abstract.} Recently, the geodesibility of planar vector fields, which are
algebrizable (differentiable in the sense of Lorch for some
associative and commutative unital algebra), has been established.
In this paper, we consider algebrizable three-dimensional vector
fields, for which we give rectifications and Riemannian metrics
under which they are geodesible. Furthermore, for each of these
vector fields $F$ we give two first integrals $h_1$ and $h_2$ such
that the integral curves of $F$ are locally defined by the
intersections of the level surfaces of $h_1$ and $h_2$.
\end{small}

\noindent \textbf{Keyword}: \texttt{Vector fields},
\texttt{Riemannian metrics}, \texttt{Lorch differentiability},
\texttt{Geodesible vector fields},



\section*{Introduction}

The Cauchy integral theorem is satisfied on algebras, for a proof
by using analytic arguments, see \cite{Ket} and \cite{Lor};
\cite{EL1} and \cite{PhiAlgDE2018} for proofs by using
conservative vector fields. Thus, the fundamental theorem of
calculus on algebras is satisfied. The theory of analytic
functions over algebras extends other results of the classical
complex function theory; see \cite{Blum}, \cite{Ket}, \cite{Lor},
\cite{Shef}, \cite{War1}, and \cite{War2}. Differentiability with
respect to finite dimensional associative algebras over $\mathbb
R$ (not necessarily commutative) has been recently reconsidered in
\cite{JSC} under the name of $\mathcal A$-Calculus.

In this paper, for fluency in writing planar vector fields,
functions with domain and values in $\mathbb R^n$ for $n=2,3$ are
denoted in the same way.

Let $\mathbb A$ be an \textit{algebra}; this means $\mathbb A$ is
the $\mathbb R$-linear space $\mathbb R^n$ endowed with a
structure of an associative, commutative algebra with unit that is
denoted by $e$; see \cite{Pie}. An element $a\in\mathbb A$ is
\textit{regular} when its multiplicative inverse $a^{-1}$ exists,
which is also denoted by $\frac{e}{a}$. The set of all the regular
elements in $\mathbb A$ is denoted by $\mathbb A^*$.

We say that a vector field $F:\Omega\subset\mathbb R^n\to\mathbb
R^n$ is $\mathbb A$-\textit{algebrizable} on an open set $\Omega$
when $F$ is differentiable in the sense of Lorch with respect to
$\mathbb A$ on $\Omega$. We simply say $F$ is
\textit{algebrizable} if $F$ is $\mathbb A$-algebrizable for some
algebra $\mathbb A$. The Lorch differential of $F$ with respect to
$\mathbb A$ at a point $w$ is called $\mathbb A$-\emph{derivative}
of $F$ at $w$, and it is denoted by $F'(w)$. Thus, if $F$ is
$\mathbb A$-algebrizable on $\Omega$, then $F'$ is the function
defined by the $\mathbb A$-derivative of $F$. This definition is
the same for each function $H:\Omega\subset\mathbb R^n\to\mathbb
R^n$. An $\mathbb A$-\emph{antiderivative} $H$ of a vector field
$F$ is an $\mathbb A$-algebrizable function $H$ such that $H'=F$.

In this work  $F=(f_1,\dots,f_n)$ denotes a vector field defined
on an open set $\Omega\subset\mathbb R^n$ for $n=2,3$ and $\mathbb
A$ an algebra such that $F$ is $\mathbb A$-algebrizable. We call
the \emph{regular domain} of $F$ with respect to $\mathbb A$ to
the set $\Omega(F,\mathbb A^*):= F^{-1}(\mathbb A^*)$, which is
characterized in Subsection \ref{rdof}. Thus, $F(\mathbb
A\setminus\mathbb A^*)$ is the set of the singular points of $F$
with respect to $\mathbb A$.

In \cite{AMEE} conditions are given under which two-dimensional
vector fields are algebrizable. In this paper we rewrite these
conditions in terms of the generalized Cauchy-Riemann equations.
Theorem \ref{A0} and Proposition \ref{algebra3pl} describe
families of two-dimensional algebras. Each algebrizable planar
vector field is algebrizable for some of these algebras. Thus, in
order to consider all the algebrizable planar vector fields we do
not have to classify the two dimensional algebras. If we consider
a system of generalized Cauchy-Riemann equations associated with
an arbitrary two-dimensional algebra $\mathbb A$, by making
elementary operations with equations it fall in some of the
families of generalized Cauchy-Riemann equations considered in
Lemmas \ref{ecr12} and \ref{ecr13}. Thus, if a given vector field
is algebrizable with respect to $\mathbb A$, then it must be
algebrizable with respect to an algebra in the families considered
in this paper.

Conditions are given in \cite{AMEE3} under which three-dimensional
vector fields are algebrizable. We rewrite these results in terms
of generalized Cauchy-Riemann equations, see Theorems \ref{t1},
\ref{t2}, and \ref{t3}. The systems of six generalized
Cauchy-Riemann equations satisfying conditions given in Theorem 2
of \cite{War2} are linearly equivalent to some of the systems
given in Lemma \ref{ecr123}. Thus, the vector fields satisfying
some of the systems of generalized Cauchy-Riemann equations given
in \cite{War2} are algebrizable with respect to some of the
algebras given in Theorem \ref{A1}. The three-dimensional algebras
considered in this paper are given in Theorem \ref{A1}, and
Propositions \ref{eas}, \ref{eas3}.

The Cauchy integral theorem is fundamental in this work, since a
line integral with respect to a three dimensional algebra $\mathbb
A$ of $\frac{e}{F}$ for a three-dimensional vector field $F$, can
be written by
\begin{equation}\label{g12}
\int_\gamma\frac{e}{F}d\xi=\left(\int_\gamma G_1\cdot
ds,\int_\gamma G_2\cdot ds,\int_\gamma G_3\cdot ds\right),
\end{equation}
for every differentiable curve $\gamma$ contained in
$\Omega(F,\mathbb A^*)$, as well as for certain conservative
vector fields $G_1$, $G_2$, $G_3$; see Propositions \ref{palvf},
\ref{proprst}, \ref{ancbw}, and \ref{aecu123}. The existence of
$\mathbb A$-antiderivatives of $\mathbb A$-algebrizable functions
and vector fields is guaranteed by the fundamental theorem of
calculus on algebras: each $\mathbb A$-algebrizable vector field
has an $\mathbb A$-antiderivative on each simply connected region.

It is said that a vector field $F$ on $\Omega$ is
\textit{geodesible} on $\Omega(F,\mathbb A^*)$ if there exists a
Riemannian metric on $\Omega(F,\mathbb A^*)$ for which the
integral curves of $F$ are unitary geodesics. The geodesibility of
algebrizable planar vector fields was established in \cite{Dyg}.

The function $H:\Omega(F,\mathbb A^*)\to\mathbb R^3$ defined
locally by the line integral
\begin{equation}\label{isolo}
H(w)=\int_\gamma \frac{e}{F}d\xi,
\end{equation}
is locally an $\mathbb A$-antiderivative of $\frac{e}{F}$, where
$\gamma$ is a differentiable curve in $\Omega(F,\mathbb A^*)$
starting at $w_0$ and ending at $w$. In this paper, we introduce
$H$ for $\mathbb A$-algebrizable three-dimensional vector fields
$F$. The function $H$ could be multivalued if the homotopy group
of some of the connected components of $\Omega(F,\mathbb A^*)$ is
not trivial. If all the connected components of $\Omega(F,\mathbb
A^*)$ are simply connected, then $H$ is a univalued function
defined on $\Omega(F,\mathbb A^*)$. In addition, it depends only
on $F$ and a base point $w_0$ in every connected component of
$\Omega(F,\mathbb A^*)$. For a point $w_0\in\Omega(F,\mathbb A^*)$
and for each of the multiple values of $H(w_0)$, the function $H$
can be restricted to a local neighborhood of $w_0$ such that the
restricted function $H_{\alpha}$ will be a diffeomorphism and
$H_{\alpha*}F(w)=e$ for all $w$ in the domain of $H_{\alpha}$,
where $e$ is the unit of $\mathbb A$; see Corollary \ref{ecalpra}.

Denote by $\lambda$ the natural Riemannian metric of $\mathbb
R^3$. The pullback of $\lambda$ under $H$ associates with $F$ the
Riemannian metric $g=\frac{1}{\|e\|^2}H^*\lambda$ defined on
$\Omega(F,\mathbb A^*)$, where $\|e\|$ denotes the Euclidean norm
of the unit $e$ of $\mathbb A$. In the rest of the paper $g$ will
denote this Riemannian metric.

In this paper we consider the vector fields $E_i=e_iF$ for
$i=1,2,3$. We show that $\langle E_i,G_j\rangle=\delta_{ij}$ for
$i,j\in\{1,2,3\}$ with respect to the usual inner product of
$\mathbb R^3$, where $\delta_{ij}$ denotes the Kronecker delta.
Thus, first integrals $h_i$ for $E_j$ for $i,j\in\{1,2,3\}$ and
$i\neq j$ are the components of an antiderivative
$H=(h_1,h_2,h_3)$ of $\frac{e}{F}$, which is given in (\ref{g12})
and (\ref{isolo}). We give some relations between $F$, $E_i$,
$G_i$, $h_i$, for $i=1,2,3$; see Propositions \ref{proprst},
\ref{ancbw}, and \ref{aecu123}.

We say that two differentiable functions
$h_1,h_2:\Omega\subset\mathbb R^3\to\mathbb R$ defined on an open
set $\Omega\subset\mathbb R^3$ have \emph{transversal
intersections}, if for every pair of level surfaces
$S_1=h_1^{-1}(c_1)$, $S_2=h_2^{-1}(c_2)$ with $S_1\cap
S_2\neq\emptyset$, $S_1$ and $S_2$ have transversal intersection.
Each algebrizable three-dimensional vector field $F$ defined on an
open set $\Omega$ admits two first integrals $h_1$, $h_2$ defined
locally on $\Omega(F,\mathbb A^*)$ that have transversal
intersections; see Proposition \ref{calprin}. The image of the
integral curves of $F$ are locally contained in the transversal
intersections of the level surfaces of $h_1$ and $h_2$; see
Corollary \ref{titfi}.

Since we are working in an algebraic framework, antiderivatives of
$\frac{e}{F}$ can be calculated by using the integration rules
known for functions of a single real or complex variable. For
example, if $\frac{e}{F(w)}=w^{-n}$, then an antiderivative of
$\frac{e}{F}$ is $H(w)=\frac{-e}{(n-1)w^{n-1}}$; see Example
\ref{pipp}. Thus, for a family of vector fields $F$, their
antiderivatives $H$ can be calculated by making operations in
their corresponding algebras.

In this paper we define a family of three-dimensional vector
fields
$$
\mathcal F=\{F_{\theta,\phi}\,:\,0\leq\theta\leq
2\pi,\,0\leq\phi\leq \pi\},
$$
where $F_{\theta,\phi}$ is the product
$F_{\theta,\phi}=u_{\theta,\phi}F$ of
$u_{\theta,\phi}=(\cos\theta\sin\phi,\sin\theta\sin\phi,\cos\phi)$
and $F$ with respect to $\mathbb A$. Thus,
$E_1=F_{0,\frac{\pi}{2}}$, $E_2=F_{\frac{\pi}{2},\frac{\pi}{2}}$,
and $E_3={F_{0,0}}$. All the vector fields on $\mathcal F$ are
geodesible with respect to $g$, and they are commutative: the Lie
bracket of each pair of vector fields on $\mathcal F$ is equal to
zero. By the definition of $g$, we have that $\{E_1,E_2,E_3\}$ is
orthonormal with respect to $g$, $g$ is flat, and $g$ can be
written in terms of the components of $F$; see Proposition
\ref{tenmr}.

If $N\subset\Omega(F,\mathbb A^*)$ is a simply connected
neighborhood and it is small enough, then the metric $d_F:N\times
N\to\mathbb R$ arising from $g$ satisfies
$$
d_F(w,w')=\|H(w)-H(w')\|,
$$
where $\|\cdot\|$ denotes the Euclidean norm in $\mathbb R^2$; see
Theorem \ref{mh}. We give two examples in Subsection \ref{med}.

H. Gluck \cite{Gluk} and D. Sullivan \cite{Sull} recognized the
obstructions of a non-vanishing vector field $F$ to be geodesible.
In the case of complex analytic vector fields $F$, the theory of
quadratic differentials (as in J. A. Jenkins \cite{Jen}, K.
Strebel \cite{S3}, J. Muciño-Raymundo and C. Valero-Valdéz
\cite{Mur2}) recognizes that $H$ given by (\ref{isolo}) is a local
rectification. In general, for each algebrizable planar vector
field $H$ is a multivalued rectification of $F$, when restricted
to each simply connected domain defines a local rectification for
each of the multiple values of $H$; see \cite{EL1}.

In Section \ref{s1}, we consider $n$-dimensional algebras for
$n=2,3$ and their generalized Cauchy-Riemann equations. In Section
\ref{s2}, the algebrizability of $n$-dimensional vector fields for
$n=2,3$ is characterized in terms of the generalized
Cauchy-Riemann equations. In Section \ref{s3}, we recall the usual
line integral, define a line integral on algebras, and give
multivalued rectifications of algebrizable vector fields. In
Section \ref{camposfygs}, for each algebrizable vector field $F$
we give the expressions of the vector fields $E_i$ and $G_i$ for
$i=1,2,3$. In addition, first integrals $h_j$, $h_k$ of $G_i$ for
$\{i,j,k\}=\{1,2,3\}$ are given. We show that each antiderivative
$H$ of $F$ gives first integrals of $F$ and that the integral
curves of $F$ are locally contained in the intersection of the
level surfaces of $h_j$ and $h_k$. In Section \ref{rmcof}, it is
showed that $\{E_1,E_2,E_3\}$ is an orthonormal frame for the
Riemannian metric $g$, and the components of $g$ are expressed in
terms of the solutions of a linear system of six equations with
six variables. Furthermore, the metric $d_F$ locally arising from
$g$ is expressed in terms of $H$.

\section{Algebras and algebrizability}\label{s1}

\subsection{First fundamental representations of algebras}\label{s11}

Consider an algebra $\mathbb A$ and a base
$\beta=\{\beta_1,\cdots,\beta_n\}$ of $\mathbb A$. The product
between the elements of $\beta$ is given by
$\beta_i\beta_j=\sum_{k=1}^n c_{ijk}\beta_k$ where
$c_{ijk}\in\mathbb R$ for $i,j,k\in\{1,\cdots,n\}$ are called
\emph{structure constants of $\mathbb A$ associated with $\beta$}.
The injective linear homomorphism $R_\beta:\mathbb A\to
M(n,\mathbb R)$ defined by $R_\beta(\beta_i)=R_i$, where $R_i$ is
the matrix with $[R_i]_{jk}=c_{ikj}$, for $i=1,2,\cdots,n$, is
called the \emph{first fundamental representation of $\mathbb A$
associated with $\beta$}. If $\mathbb A$ as a linear space is
$\mathbb R^n$, $R:\mathbb A\to M(n,\mathbb R)$ will denote the
first fundamental representation of $\mathbb A$ associated with
the ordered standard basis $\beta=\{e_1,\cdots,e_n\}$, and in this
case $R$ will be called simply a \emph{first fundamental
representation of $\mathbb A$}. When we refer to an algebra
$\mathbb A$, we are considering the linear space $\mathbb R^n$
with an algebra structure.

\subsection{Two-dimensional algebras}\label{aep}

In this subsection, we give quit a different presentation of the
algebras given in \cite{AMEE} which are used in \cite{Dyg}.

\bthm\label{A0} Let $p_1$ and $p_2$ be parameters and
$\{e_r,e_s\}=\{e_1,e_2\}$. The linear space $\mathbb{R}^{2}$
endowed with the product
\begin{equation}\label{algebra0}
  \begin{tabular}{c|cc}
  $\cdot$ & $e_r$ & $e_s$ \\
  \hline
  $e_r$ & $e_r$ & $e_s$ \\
  $e_s$ & $e_s$ & $p_1e_r+p_2e_s$ \\
\end{tabular}
\end{equation}
is an algebra $\mathbb{A}$ with unit $e=e_r$ and with structure
constants
$$
\begin{array}{cccc}
  c_{rrr}=1, & c_{rsr}=0, & c_{srr}=0, & c_{ssr}=p_1, \\
  c_{rrs}=0, & c_{rss}=1, & c_{srs}=1, & c_{sss}=p_2. \\
\end{array}%
$$
\ethm \noindent\textbf{Proof.} Consider the matrix
$$
I=\left(%
\begin{array}{cc}
  1 & 0 \\
  0 & 1 \\
\end{array}%
\right),\quad A=\left(%
\begin{array}{cc}
  0 & p_1 \\
  1 & p_2 \\
\end{array}%
\right),\quad \quad B=\left(%
\begin{array}{cc}
  p_2 & 0 \\
  p_1 & 1 \\
\end{array}%
\right).
$$
If $s=2$ we take $M_r=I$, $M_s=A$ and $M_s=B$, $M_r=I$ in the
other case. Therefore, $M_sM_s=p_1M_r+p_2M_s$. Thus, the linear
space spanned by $\beta=\{M_1,M_2\}$ is a two-dimensional
commutative matrix algebra $\mathbb M$ with a first fundamental
representation $R_\beta$ that is the identity automorphism.
Consider the linear space $\mathbb R^2$ endowed with the product
(\ref{algebra0}) and $R:\mathbb R^2\to\mathbb M$ be defined by
$R(x_re_r+x_se_s)=x_rM_r+x_sM_s$, then $R$ is a linear
isomorphism. So, $\mathbb R^2$ is endowed with the algebra
structure of $\mathbb M$ induced by $R$. Thus, we get an algebra
$\mathbb A$ with a product given by (\ref{algebra0}), in which the
first the fundamental representation is given by $R$. $\Box$

We use the notation $\mathbb A^2_r(p_1,p_2)$ for the algebra with
unit $e=e_r$ and parameters $p_1,p_2$, as defined in Theorem
\ref{A0}.

\bprop\label{algebra3pl} The linear space $\mathbb R^2$ endowed
with the product $e_ie_i=e_i$ for $i=1,2$, and $e_ie_j=0$ in the
other case, is an algebra $\mathbb A$ with unit $e=e_1+e_2$.
\eprop

We use the notation $\mathbb A^2_{1,2}$ for the algebra with unit
$e=e_1+e_2$, as defined in Proposition \ref{algebra3pl}.

\subsection{Three-dimensional algebras}

In this subsection, we give quit a different presentation of the
algebras given in \cite{AMEE3}.

\bthm\label{A1} Let $p_1,p_2,\cdots,p_9$ be constants satisfying
the equalities
\begin{equation}\label{ceq}
    \begin{array}{c}
      p_7 = p_1p_4+p_2p_6-p_2p_3-p_4^2, \\
      p_8= p_3p_4-p_2p_5,\qquad\qquad\,\,\,\,\,\,\\
      p_9= p_1p_5+p_3p_6-p_4p_5-p_3^2, \\
    \end{array}
\end{equation}
and $\{e_r,e_s,e_t\}=\{e_1,e_2,e_3\}$. Thus, the linear space
$\mathbb{R}^{3}$ endowed with the following product
\begin{equation}\label{algebra1}
  \begin{tabular}{c|ccc}
  $\cdot$ & $e_r$ & $e_s$ & $e_t$ \\
  \hline
  $e_r$ & $e_r$ & $e_s$ & $e_t$ \\
  $e_s$ & $e_s$ & $p_7e_r+p_1e_s+p_2e_t$ & $p_8e_r+p_3e_s+p_4e_t$ \\
  $e_t$ & $e_t$ & $p_8e_r+p_3e_s+p_4e_t$ & $p_9e_r+p_5e_s+p_6e_t$ \\
\end{tabular}
\end{equation}
is an algebra $\mathbb{A}$ with unit $e=e_r$. \ethm
\noindent\textbf{Proof.} Consider the matrices
$$
R_s=\left(%
\begin{array}{ccc}
  0 & p_7 & p_8 \\
  1 & p_1 & p_3 \\
  0 & p_2 & p_4 \\
\end{array}%
\right),\,\,\,R_t=\left(%
\begin{array}{ccc}
  0 & p_8 & p_9 \\
  0 & p_3 & p_5 \\
  1 & p_4 & p_6 \\
\end{array}%
\right),
$$
and $R_r=I$, where $I\in M(3,\mathbb R)$ is the identity matrix.
The following products of matrices satisfy
\begin{eqnarray*}
  R_sR_s &=& p_7R_1+p_1R_s+p_sR_t\\
  R_sR_t &=& p_8R_1+p_3R_s+p_4R_t \\
  R_tR_s &=& p_8R_1+p_3R_s+p_4R_t \\
  R_tR_t &=& p_9R_1+p_5R_s+p_6R_t.
\end{eqnarray*}
Thus, the three-dimensional space $\mathbb M$ spanned by
$\beta=\{R_1,R_2,R_3\}$ is a matrix algebra in $M(3,\mathbb R)$
with the first fundamental representation $R_\beta$ that is the
identity automorphism. Let $R:\mathbb R^3\to M(3,\mathbb R)$ be
the representation linearly defined by $R(e_i)=R_i$ for $i=1,2,3$.
We endow $\mathbb R^3$ with the algebra structure induced by $R$.
Thus, we get an algebra $\mathbb A$ with the product given in
(\ref{algebra1}), and with first fundamental representation given
by $R$. Since $R(e_r)=R_r$ is the identity matrix, we have that
$e_r$ is the unit $e$ of $\mathbb A$. $\Box$

We use the notation $\mathbb A^3_r(p_1,\cdots,p_6)$ for the
algebra with unit $e=e_r$ and parameters $p_1,\cdots,p_6$, as
defined in Theorem \ref{A1}.

We call the set of equations (\ref{ceq}) \emph{commutativity
equations}. From these, we can obtain
\begin{equation}\label{ecd}
p_3p_7+p_4p_8=p_1p_8+p_2p_9,\qquad p_5p_7+p_6p_8=p_3p_8+p_4p_9.
\end{equation}

We give more examples of algebras in the following proposition.

\bprop\label{eas} Let $\mathbb B$ be a two-dimensional algebra
with product
\begin{center}
\begin{tabular}{c|c c}
$\cdot$ & $ \beta_s$& $\beta_t$ \\ \hline $\beta_s$ & $\beta_s$ & $\beta_t$\\
 $\beta_t$ & $\beta_t$ & $ p_1\beta_s+p_2\beta_t$ \\
\end{tabular},
\end{center}
with respect to a basis $\beta=\{\beta_1,\beta_2\}$ of $\mathbb
B$, where $\{\beta_s,\beta_t\}=\{\beta_1,\beta_2\}$. If
$\{e_r,e_s,e_t\}=\{e_1,e_2,e_3\}$, then the linear space $\mathbb
R^3$ endowed with the product
\begin{center}
\begin{tabular}{c|c c c}
$\cdot$ & $e_r$ & $ e_s$& $e_t$ \\ \hline $e_r$ & $e_r$ & $ 0$& $0$ \\
$e_s$& $0$ & $e_s$ & $e_t$\\
 $e_t $&$0$ & $e_t$ & $p_1e_s+p_2e_t$ \\
\end{tabular},
\end{center}
is an algebra $\mathbb A$. \eprop \noindent\textbf{Proof.} The
matrices
$$
R_1=\left(%
\begin{array}{ccc}
  1 & 0 & 0 \\
  0 & 0 & 0 \\
  0 & 0 & 0 \\
\end{array}%
\right), \quad R_2=\left(%
\begin{array}{ccc}
  0 & 0 & 0 \\
  0 & 1 & 0 \\
  0 & 0 & 1 \\
\end{array}%
\right),\quad R_3\left(%
\begin{array}{ccc}
  0 & 0 & 0 \\
  0 & 0 & p_1 \\
  0 & 1 & p_2 \\
\end{array}%
\right)
$$
define a commutative matrix algebra $\mathbb M$ with respect to
the matrix product. If we linearly define $R:\mathbb R^3\to\mathbb
M$ by $R(e_r)=R_1$, $R(e_s)=R_2$, and $R(e_t)=R_3$, then $R$ is a
linear isomorphism that we use to define a product between the
elements of $\beta=\{e_r,e_s,e_t\}$ in such a way that $\mathbb
R^3$ is an algebra $\mathbb A$ with first fundamental
representation $R_\beta=R$. Thus, the proof is finished. $\Box$

We use the notation $\mathbb A^3_{r,s}(p_1,p_2)$ for the algebra
with unit $e=e_r+e_s$ and parameters $p_1,p_2$, as defined in
Proposition \ref{eas}.

The last type of three-dimensional algebra considered in this
paper is the given in the following proposition.

\bprop\label{eas3} The linear space $\mathbb R^3$ endowed with the
product $e_ie_i=e_i$ for $i=1,2,3$, and $e_ie_j=0$ in the other
cases, is an algebra $\mathbb A$. \eprop

We use the notation $\mathbb A^3_{1,2,3}$ for the algebra with
unit $e=e_1+e_2+e_3$, as defined in Proposition \ref{eas3}.

\subsection{Algebrizability}\label{s12}

The definitions and results regarding algebrizability, classical
line integrals and line integrals on algebras, of vector fields
and functions $F:\Omega\subset\mathbb R^n\to\mathbb R^n$ are the
same in this paper. These vector fields and functions are denoted
in the same way $F=(f_1,\cdots,f_n)$.

Now, we recall the definition of algebrizability of vector fields.

\bdefe\label{55} Let $\mathbb A$ be an algebra and
$F:\Omega\subset\mathbb R^n\to\mathbb R^n$ a vector field defined
on an open set $\Omega$. We say $F$ is $\mathbb
A$-\textit{algebrizable} on $\Omega$ if there exists a function
$F':\Omega\subset\mathbb R^n \rightarrow\mathbb R^n$, called
\emph{$\mathbb A$-derivative} of $F$ on $\Omega$, which satisfies
\begin{equation}\label{difa}
\lim_{\xi\to 0}\frac{|F(w_0+\xi)-F(w_0)-F'(w_0)\xi|}{|\xi|}=0
\end{equation}
for all $w_0\in\Omega$, where $F'(w_0)\xi$ denotes the product in
$\mathbb A$ of $F'(w_0)$ and $\xi$. \edefe

We simply say $F$ is \emph{algebrizable} if $F$ is $\mathbb
A$-algebrizable for some algebra $\mathbb A$.

If $F:\Omega\subset\mathbb R^n\rightarrow\mathbb R^n$ is
differentiable in the usual sense, we denote by $JF:\Omega\to
M(n,\mathbb R)$ the function defined by the Jacobian matrix of $F$
in $w$: the $(i,j)$ entry of $JF(w)$ is
$[JF(w)]_{i,j}=\frac{\partial f_i}{\partial x_j}|_{(w)}$.

The following two results are given in papers related to
differentiability on algebras.

\blem\label{l1} Let $\mathbb A$ be an algebra with the first
fundamental representation $R:\mathbb R^n\to M(n,\mathbb R)$. For
each $u\in\mathbb A$, we have that $uw=[R(u)w^T]^T$, for all
$w\in\mathbb A$, where $uw$ denotes the product in $\mathbb A$,
$R(u)w^T$ denotes the matrix product of $R(u)$ by $w^T$, and $^T$
denotes the transpose matrix. \elem \noindent\textbf{Proof.} It is
sufficient to show that $e_ie_j=[R(e_i)e_j^T]^T$, which can be
shown in a straightforward form. $\Box$

\blem\label{ljd} (\cite{War2}) Let $F:\Omega\subset\mathbb R^n
\to\mathbb R^n$ be a vector field, where $\Omega$ is an open set
and $\mathbb A$ is an algebra. Therefore, $F$ is $\mathbb
A$-algebrizable on $\Omega$ if and only if $F$ is differentiable
in the sense of Fréchet on $\Omega$ and the function $JF:\Omega\to
M(n,\mathbb R)$ has an image in $R(\mathbb A)$, where $R$ is the
first fundamental representation of $\mathbb A$. \elem
\noindent\textbf{Proof.} Suppose that $F$ is $\mathbb
A$-algebrizable at $w_0$. The $\mathbb A$-derivative and the usual
differential of $F$ at $w_0$ satisfy the equality
$F'(w_0)w=(JF(w_0)w^T)^T$, where $F'(w_0)w$ denotes the product on
$\mathbb A$ of $F'(w_0)$ by $w$, $JF(w_0)w^T$ is the matrix
product of $JF(w_0)$ by $w^T$, and $^T$ denotes the transpose
matrix. From Lemma \ref{l1}, we have that $JF(w_0)=R(F'(w_0))$;
thus $JF(w_0)\in R(\mathbb A)$.

Suppose that $JF(w_0)\in R(\mathbb A)$, then by Lemma \ref{l1}
$v=R^{-1}(JF(w_0))$ satisfies $vw=(JF(w_0)w^T)^T$. Since $F$ is
differentiable in the usual sense, we have
$$
0=\lim_{\xi\to
0}\frac{|F(w_0+\xi)-F(w_0)-[JF(w_0)\xi^T]^T|}{|\xi|}=\lim_{\xi\to
0}\frac{|F(w_0+\xi)-F(w_0)-v\xi|}{|\xi|}:
$$
the $\mathbb A$-derivative of $F$ at $w_0$ is $F'(w_0)=v$. $\Box$

\subsection{Generalized Cauchy-Riemann equations}\label{s13}

In this subsection, we give a characterization of the
algebrizability based on the generalized Cauchy-Riemann equations.

Given an algebra $\mathbb A$ with structure constants $c_{ijk}$,
the set of partial differential equations
\begin{equation}\label{cre}
\sum_{i=1}^n c_{ijk}\frac{\partial f_i}{\partial x_m}=\sum_{i=1}^n
c_{imk}\frac{\partial f_i}{\partial x_j},\quad k=1,2,\cdots,n,
\end{equation}
appears in \cite{Ket} p. 646 and is called \emph{generalized
Cauchy-Riemann equations} of $\mathbb A$. A vector field
$F=(f_1,\cdots,f_n)$ is $\mathbb A$-algebrizable on $\Omega$ if
and only if $\{f_1,\cdots,f_n\}$ satisfy (\ref{cre}) on $\Omega$.
We use GCRE for generalized Cauchy-Riemann equations.

\bdefe We say that two linear systems of partial differential
equations are \textbf{linearly equivalent} if under a finite
number of elementary operations we can obtain one system from the
other. \edefe

We think that each three-dimensional algebra $\mathbb A$ has GCRE
that are linearly equivalent to the GCRE of some of the
three-dimensional algebras given in this section however, we do
not prove them.

\section{On algebrizability of vector fields}\label{s2}

We use the notations for algebras introduced in Section \ref{s1}.

\subsection{Algebrizability of planar vector fields}

In this subsection, we consider planar vector fields $F=(f_1,f_2)$
defined in regions $\Omega$ of $\mathbb R^2$.

In the following lemma, we give the GCRE associated with $\mathbb
A^2_r(p_1,p_2)$.

\blem\label{ecr12} Let $\mathbb A=\mathbb A^2_r(p_1,p_2)$ and
$\{r,s\}=\{1,2\}$. The GCRE of $\mathbb A$ are given by
\begin{equation}\label{ECR2}
    \frac{\partial f_r}{\partial x_s} = p_1\frac{\partial f_s}{\partial
    x_r},\qquad \frac{\partial f_s}{\partial x_s} = \frac{\partial f_r}{\partial
  x_r}+p_2\frac{\partial f_s}{\partial x_r}.
\end{equation}
\elem \noindent\textbf{Proof.} Let $\mathbb A=\mathbb
A^2_r(p_1,p_2)$. The GCRE (\ref{cre}) with respect to the
structure constants of $\mathbb A$ given in Theorem \ref{A0}
define the set of equations (\ref{ECR2}).$\Box$

The corollary of Theorem 2 of \cite{War2} considers linear systems
of two partial differential equations and gives conditions in
order to define algebras. Each system satisfying these conditions
is linearly equivalent to some of those considered in Lemma
\ref{ecr12}. Thus, if some vector field $F=(f_1,f_2)$ satisfies a
linear system of partial differential equations satisfying the
conditions of Theorem 2 of \cite{War2}, then $F$ satisfies some of
the systems of GCRE given in Lemma \ref{ecr12}.

If $\mathbb A=\mathbb A^2_r(p_1,p_2)$, the $\mathbb
A$-algebrizability of a planar vector field $F$ can be verified by
the GCRE (\ref{ECR2}), as we see in the following theorem.

\bthm\label{t10} Let $F:\Omega\subset\mathbb R^{2}\to \mathbb
R^{2}$ be differentiable in the sense of Fréchet over an open set
$\Omega$ and $\mathbb A=\mathbb A^2_r(p_1,p_2)$. Therefore, the
components of $F$ satisfy equations (\ref{ECR2}) if and only if
$F$ is $\mathbb{A}$-algebrizable. \ethm \noindent\textbf{Proof.}
Let $\mathbb A=\mathbb A^2_r(p_1,p_2)$. Suppose that the
components of $F$ satisfy equations (\ref{ECR2}), then
$$
Jf=\left(%
\begin{array}{cc}
  1 & 0 \\
  0 & 1 \\
\end{array}%
\right)\frac{\partial f_1}{\partial x_1}+\left(%
\begin{array}{cc}
  0 & p_1 \\
  1 & p_2 \\
\end{array}%
\right)\frac{\partial f_2}{\partial x_1}.
$$
By proof of Theorem \ref{A0}, we have that $JF\in R(\mathbb A)$.
Thus, by Lemma \ref{ljd} we have that $F$ is $\mathbb
A$-algebrizable.

Suppose that $F$ is $\mathbb A$-algebrizable, then by proof of
Theorem \ref{A0} and by Lemma \ref{ljd} we have that
$$
JF=\left(%
\begin{array}{cc}
  1 & 0 \\
  0 & 1 \\
\end{array}%
\right)u_1+\left(%
\begin{array}{cc}
  0 & p_1 \\
  1 & p_2 \\
\end{array}%
\right)u_2
$$
for certain functions $u_1$ and $u_2$ depending on $(x_1,x_2)$.
The last equality for $JF$ implies that $u_i=\frac{\partial
f_i}{\partial x_1}$ for $i=1,2$. So, the components of $F$ satisfy
equations (\ref{ECR2}).

The other case can be proved in a similar way. $\Box$

In the following lemma, we give the GCRE associated with $\mathbb
A^2_{1,2}$ .

\blem\label{ecr13} The GCRE of $\mathbb A=\mathbb A^2_{1,2}$ are
given by
\begin{equation}\label{CRE2d}
    \frac{\partial f_1}{\partial x_2} = 0,\qquad \frac{\partial f_2}{\partial x_1} =0.
\end{equation}
\elem \noindent\textbf{Proof.} The structure constants of $\mathbb
A^2_{1,2}$ are $c_{iii}=1$ for $i=1,2$ and $c_{ijk}=0$ in another
case. By substituting these constants in equations (\ref{cre}), we
obtain (\ref{CRE2d}). Thus, the proof is finished. $\Box$

The set of GCRE (\ref{CRE2d}) does not satisfy conditions of the
corollary of Theorem 2 of \cite{War2}. In \cite{AMEE} it is shown
with quit a different presentation of the subject that every
algebrizable planar vector field $F=(f_1,f_2)$ satisfies equations
(\ref{ECR2}) for a suitable choice of parameters $p_1,p_2$ and
$\{r,s\}=\{1,2\}$, or it satisfies equations (\ref{CRE2d}). This
gives a criterion that serves to determine the algebrizablity of a
planar vector field $F$.

If $\mathbb A=\mathbb A^2_{1,2}$, the $\mathbb A$-algebrizability
of a planar vector field $F$ can be verified by the GCRE
(\ref{CRE2d}), as we see in the following theorem.

\bthm\label{t11} Let $F:\Omega\subset\mathbb R^{2}\to \mathbb
R^{2}$ be differentiable in the sense of Fréchet over an open set
$\Omega$ and $\mathbb A=\mathbb A^2_{1,2}$. Thus, the components
of $F$ satisfy equations (\ref{CRE2d}) if and only if $F$ is
$\mathbb{A}$-algebrizable. \ethm \noindent\textbf{Proof.} The
proof is similar to that of Theorem \ref{t10}.$\Box$

For $F$ $\mathbb A$-differentiable with respect to $\mathbb
A=\mathbb A^2_{1,2}$, the component $f_r$ is a differentiable
function of the variable $x_r$, and it does not depend on $x_s$
for $\{r,s\}=\{1,2\}$.

\subsection{Algebrizability of three-dimensional vector fields}

In this subsection we study the algebrizability of vector fields
$F$ defined on regions $\Omega$ of $\mathbb R^3$. We characterize
the algebrizability of families of three-dimensional vector fields
in terms of GCRE. In the rest of the paper, $F$ will denote a
three-dimensional vector field $F=(f_1,f_2,f_3)$.

Now consider $F:\Omega \subset \mathbb{R}^3 \rightarrow
\mathbb{R}^3$ being differentiable in the sense of Fréchet on an
open set $\Omega$. Writing $F$ in terms of its components, the
usual matrix-valuated function $JF:\Omega \subset \mathbb{R}^3
\rightarrow M(3,\mathbb R)$ defined by the partial derivatives is
given by $[JF]_{ij}=\frac{\partial f_i}{\partial x_j}$.

\blem\label{ecr123} Let $\mathbb A$ be the algebra $\mathbb
A^3_r(p_1,\cdots,p_6)$ and $\{e_r,e_s,e_t\}=\{e_1,e_2,e_3\}$. The
GCRE of $\mathbb A$ are given by
\begin{equation}\label{ECR123}
  \begin{array}{cc}
    \frac{\partial f_r}{\partial x_s}=p_7\frac{\partial f_s}{\partial x_r}+p_8\frac{\partial f_t}{\partial x_r},\quad\,\,\,\,\, & \frac{\partial f_r}{\partial x_t}=p_8\frac{\partial f_s}{\partial x_r}+p_9\frac{\partial f_t}{\partial x_r}, \\
    \frac{\partial f_s}{\partial x_s}=\frac{\partial f_r}{\partial x_r}+p_1\frac{\partial f_s}{\partial x_r}+p_3\frac{\partial f_t}{\partial x_r}, & \frac{\partial f_s}{\partial x_t}=p_3\frac{\partial f_s}{\partial x_r}+p_5\frac{\partial f_t}{\partial x_r}, \\
    \frac{\partial f_t}{\partial x_s}=p_2\frac{\partial f_s}{\partial x_r}+p_4\frac{\partial f_t}{\partial x_r},\quad\,\,\,\, & \quad\,\,\,\,\frac{\partial f_t}{\partial x_t}=\frac{\partial f_r}{\partial x_r}+p_4\frac{\partial f_s}{\partial x_r}+p_6\frac{\partial f_t}{\partial x_r}, \\
  \end{array}
\end{equation}
where $p_7$, $p_8$, and $p_9$ are defined by the commutativity
equations (\ref{ceq}). \elem \noindent\textbf{Proof.} Consider the
algebra $\mathbb A^3_r(p_1,\cdots,p_6)$ with unit $e_r$ for
$r\in\{1,2,3\}$. The six equations (\ref{cre}) obtained with $j=r$
and $m\in\{1,2,3\}\setminus\{r\}$ are the equations given above.
$\Box$

In Theorem 2 of \cite{War2} are given sufficient conditions in
order to ensure that a given linear systems of partial
differential equations is a set of GCRE for some algebra. Each
system satisfying these conditions is linearly equivalent to some
system like (\ref{ECR123}). Thus, if some vector field
$F=(f_1,f_2,f_3)$ satisfies a linear systems of partial
differential equations satisfying the conditions of Theorem 2 of
\cite{War2}, then $F$ satisfies some of the systems like
(\ref{ECR123}).

For algebras $\mathbb A=\mathbb A^3_r(p_1,\cdots,p_6)$, the GCRE
associated with $\mathbb A$ serve as a criteria for determining
the $\mathbb A$-algebrizability.

\bthm\label{t1} Let $F:\Omega\subset\mathbb R^{3}\to \mathbb
R^{3}$ be differentiable in the sense of Fréchet over an open set
$\Omega$ and $\mathbb A=\mathbb A^3_r(p_1,\cdots,p_6)$.  Thus, the
components of $F$ satisfy equations (\ref{ECR123}) if and only if
$F$ is $\mathbb{A}$-algebrizable. \ethm \noindent\textbf{Proof.}
Let $\mathbb A=\mathbb A^3_r(p_1,\cdots,p_6)$. Suppose that the
components of $F$ satisfy equations (\ref{ECR123}), then
$$
JF=\left(%
\begin{array}{ccc}
  1 & 0 & 0 \\
  0 & 1 & 0 \\
  0 & 0 & 1 \\
\end{array}%
\right)\frac{\partial f_1}{\partial x_1}+\left(%
\begin{array}{ccc}
  0 & p_7 & p_8 \\
  1 & p_1 & p_3 \\
  0 & p_2 & p_4 \\
\end{array}%
\right)\frac{\partial f_2}{\partial x_1}
+\left(%
\begin{array}{ccc}
  0 & p_8 & p_9 \\
  0 & p_3 & p_5 \\
  1 & p_4 & p_6 \\
\end{array}%
\right)\frac{\partial f_3}{\partial x_1}.
$$
By proof of Theorem \ref{A1}, we have that $JF\in R(\mathbb A)$.
Thus, by Lemma \ref{ljd} we have that $F$ is $\mathbb
A$-algebrizable.

Suppose that $F$ is $\mathbb A$-algebrizable, then by proof of
Theorem \ref{A1} and by Lemma \ref{ljd} we have that
$$
JF=\left(%
\begin{array}{ccc}
  1 & 0 & 0 \\
  0 & 1 & 0 \\
  0 & 0 & 1 \\
\end{array}%
\right)u_1+\left(%
\begin{array}{ccc}
  0 & p_7 & p_8 \\
  1 & p_1 & p_3 \\
  0 & p_2 & p_4 \\
\end{array}%
\right)u_2
+\left(%
\begin{array}{ccc}
  0 & p_8 & p_9 \\
  0 & p_3 & p_5 \\
  1 & p_4 & p_6 \\
\end{array}%
\right)u_3
$$
for certain functions $u_1,u_2,u_3$ depending on $(x_1,x_2,x_3)$.
The last equality for $JF$ implies that $u_i=\frac{\partial
f_i}{\partial x_1}$ for $i=1,2,3$. Thus, the components of $F$
satisfy equations (\ref{ECR123}).

The remaining cases can be proved in a similar way. $\Box$

\blem\label{a3dd} The GCRE of $\mathbb A^3_{r,s}(p_1,p_2)$ are
given by
\begin{equation}\label{A3DD}
\begin{array}{ccc}
  \frac{\partial f_r}{\partial x_s}=0, & \frac{\partial f_s}{\partial x_r}=0, & \frac{\partial f_s}{\partial x_t}
  =p_1\frac{\partial f_t}{\partial x_s},\qquad \\
  \frac{\partial f_r}{\partial x_t}=0, & \frac{\partial f_t}{\partial x_r}=0, & \frac{\partial f_t}{\partial x_t}
  =\frac{\partial f_s}{\partial x_s}+p_2\frac{\partial f_t}{\partial x_s}. \\
\end{array}%
\end{equation}
\elem \noindent\textbf{Proof.} We consider only the case where
$r=2$, $s=1$, and $t=3$. Then, the product of $\mathbb A$ is given
by
\begin{center}
\begin{tabular}{c|c c c}
$\cdot$ & $e_1$ & $ e_2$& $e_3$ \\ \hline $e_1$ & $e_1$ & $ 0$ & $e_3$ \\
$e_2$ & $0$ & $e_2$ & $0$\\
 $e_3$ & $e_3$ & $0$ & $p_1e_1+p_2e_3$ \\
\end{tabular},
\end{center}
from which we obtain the structure constants of $\mathbb A$:
$$
\begin{array}{ccc}
  C_{111}=1, & C_{112}=0, & C_{113}=0, \\
  C_{121}=0, & C_{122}=0, & C_{123}=0,  \\
  C_{131}=0, & C_{132}=0, & C_{133}=1,  \\
   C_{211}=0, & C_{212}=0, & C_{213}=0, \\
  C_{221}=0, & C_{222}=1, & C_{223}=0,  \\
  C_{231}=0, & C_{232}=0, & C_{233}=0,  \\
    C_{311}=0, & C_{312}=0, & C_{313}=1, \\
  C_{321}=0, & C_{322} =0, & C_{323}=0,  \\
  C_{331}=p_1, & C_{332}=0, & C_{333}=p_2.  \\
\end{array}%
$$
If we substitute these structure constants in the GCRE (\ref{cre})
of $\mathbb A$, we obtain the linear system of partial
differential equations given in the lemma for constants $r=2$,
$s=1$, and $t=3$.$\Box$

The algebras $\mathbb A^3_{r,s}(p_1,p_2)$ do not satisfy
conditions of Theorem 2 of Ward; see \cite{War2}. The GCRE
associated with $\mathbb A^3_{r,s}(p_1,p_2)$ imply
algebrizability, as we see in the following theorem.

\bthm\label{t2} Let $F:\Omega\subset\mathbb R^{3}\to \mathbb
R^{3}$ be differentiable in the sense of Fréchet over an open set
$\Omega$ and $\mathbb A=\mathbb A^3_{r,s}(p_1,p_2)$. The
components of $F$ then satisfy equations (\ref{A3DD}) if and only
if $F$ is $\mathbb{A}$-algebrizable. \ethm
\noindent\textbf{Proof.} We consider only $\mathbb
A^3_{r,s}(p_1,p_2)$ with $r=2$, $s=1$, and $t=3$, then equations
(\ref{A3DD}) are given by
$$
\begin{array}{ccc}
  \frac{\partial f_2}{\partial x_1}=0, & \frac{\partial f_1}{\partial x_2}=0, & \frac{\partial f_1}{\partial x_3}
  =p_1\frac{\partial f_3}{\partial x_1},\qquad \\
  \frac{\partial f_2}{\partial x_3}=0, & \frac{\partial f_3}{\partial x_2}=0, & \frac{\partial f_3}{\partial x_3}
  =\frac{\partial f_1}{\partial x_1}+p_2\frac{\partial f_3}{\partial x_1}. \\
\end{array}%
$$
By considering the structure constants given in Theorem
\ref{a3dd}, we have that the first fundamental representation of
$\mathbb A$ is given by
$$
R(e_1)=\left(%
\begin{array}{ccc}
  1 & 0 & 0 \\
  0 & 0 & 0 \\
  0 & 0 & 1 \\
\end{array}%
\right),\quad R(e_2)=\left(%
\begin{array}{ccc}
  0 & 0 & 0 \\
  0 & 1 & 0 \\
  0 & 0 & 0 \\
\end{array}%
\right),\quad R(e_3)=\left(%
\begin{array}{ccc}
  0 & 0 & p_1 \\
  0 & 0 & 0 \\
  1 & 0 & p_2 \\
\end{array}%
\right).
$$
We then obtain
$$
JF=\frac{\partial f_1}{\partial x_1}R(e_1)+\frac{\partial
f_3}{\partial x_1}R(e_3)+\frac{\partial f_2}{\partial x_2}R(e_2).
$$
Thus, $JF\in R(\mathbb A)$, that is, $F$ is $\mathbb
A$-algebrizable.

Suppose now that $F$ is $\mathbb A$-algebrizable, then $JF\in
R(\mathbb A)$. Thus, $JF=u_1R(e_1)+u_2R(e_2)+u_3R(e_3)$ for
certain functions $u_1$, $u_2$, and $u_3$. Therefore,
$u_2=\frac{\partial f_2}{\partial x_2}$, $f_1$, $f_3$ do not
depend on $x_2$, and $f_1$, $f_3$ satisfy equations
$$
\frac{\partial f_1}{\partial x_3}
  =p_1\frac{\partial f_3}{\partial x_1},\qquad \frac{\partial f_3}{\partial x_3}
  =\frac{\partial f_1}{\partial x_1}+p_2\frac{\partial
  f_3}{\partial
  x_1}.
$$
Thus, $f_1$, $f_2$ and $f_3$ satisfy equations (\ref{A3DD}) for
$r=2$, $s=1$, and $t=3$. $\Box$

\blem Let $\mathbb A=\mathbb A^2_{1,2}$. The GCRE associated with
$\mathbb A$ are given by
\begin{equation}\label{eccc}
\frac{\partial f_i}{\partial x_j}=0, \qquad i,j\in\{1,2,3\},\qquad
i\neq j.
\end{equation}
\elem \noindent\textbf{Proof.} The structure constants of the
algebra $\mathbb A$ given in Proposition \ref{algebra3pl} are
$c_{iii}=1$ for $i=1,2,3$ and $c_{ijk}=0$ in other case. By
substituting these constants in equations (\ref{cre}), we obtain
(\ref{eccc}). Thus, the proof is finished. $\Box$

The algebra $\mathbb A^3_{1,2,3}$ does not satisfy conditions of
Theorem 2 of Ward; see \cite{War2}. The GCRE associated with
$\mathbb A^3_{1,2,3}$ imply algebrizability, as we see in the
following theorem.

\bthm\label{t3} Let $F:\Omega\subset\mathbb R^{3}\to \mathbb
R^{3}$ be differentiable in the sense of Fréchet over an open set
$\Omega$ and $\mathbb A=\mathbb A^3_{1,2,3}$. The components of
$F$ then satisfy equations (\ref{A3DD}) if and only if $F$ is
$\mathbb{A}$-algebrizable. \ethm \noindent\textbf{Proof.} The
proof is similar to that of Theorem \ref{t1}.$\Box$

\section{On line integrals of algebrizable vector fields}\label{s3}

In the rest of the paper, we use the notation
$\langle\cdot,\cdot\rangle:\mathbb R^n\times\mathbb R^n\to\mathbb
R$ for the usual inner product, and $\mathbb A$ for an algebra
defined by the linear space $\mathbb R^3$ endowed with a product
defined in the canonical basis $\{e_1,e_2,e_3\}$ by
$e_ie_j=\sum_{k=1}^3c_{ijk}e_k$. The elements of $\mathbb R^3$ are
denoted by $w=(x_1,x_2,x_3)$, then $w=x_1e_1+x_2e_2+x_3e_3$.

\subsection{On line integrals}

Consider a vector field $F:\Omega\subset\mathbb R^3\to\mathbb R^3$
where $\Omega$ is an open set. The usual line integral of $F$ is
\begin{equation}\label{nilrn}
\int_\gamma F\cdot ds=\int_\gamma F(s)\cdot ds:=\int_{0}^a \langle
F(\gamma(t)),\dot{\gamma}(t)\rangle dt,
\end{equation}
where $\gamma:[0,a]\to\Omega$ is a differentiable curve, and
$\dot{\gamma}$ its derivative.

In the same way, the following definition of a line integral
relative to algebras is given.

\bdefe For each $\mathbb A$-algebrizable vector field $F$ defined
on an open set $\Omega$ and a differentiable curve
$\gamma:[0,a]\to\Omega$ joining $w_0$ and $w$ the \emph{line
integral of $F$ relative to $\mathbb A$} is
\begin{equation}\label{nila}
\int_\gamma F d\xi=\int_\gamma F(\xi) d\xi:=\int_{0}^a
F(\gamma(t))\dot{\gamma}(t) dt,
\end{equation}
where $F(\gamma(t))\dot{\gamma}(t)$ denotes the product of
$F(\gamma(t))$ and $\dot{\gamma}(t)$ with respect to $\mathbb A$.
\edefe

A characterization of $\Omega(F,\mathbb A^*)$ is given in
Subsection \ref{rdof}. The following proposition is satisfied for
algebrizable vector fields defined on open subsets of $\mathbb
R^n$. The proof of the first part appears in \cite{EL1} p. 198 for
the general case.

\bprop\label{palvf} Let $F$ be an $\mathbb A$-algebrizable vector
field defined on an open set $\Omega\subset\mathbb R^3$. Thus,
$$
\int_\gamma\frac{e}{F}d\xi=\left(\int_\gamma G_1\cdot
ds,\int_\gamma G_2\cdot ds,\int_\gamma G_3\cdot ds\right),
$$
for every differentiable curve $\gamma$ contained in
$\Omega(F,\mathbb A^*)$. Here $G_1$, $G_2$, and $G_3$ are
conservative vector fields (gradient) given by
$G_k=\sum_{j=1}^3\sum_{i=1}^3c_{i,j,k}g_ie_j$ for $k=1,2,3$, where
$c_{ijk}$ are the structure constants of $\mathbb A$ and $g_1$,
$g_2$, $g_3$ satisfy $(g_1,g_2,g_3)=\frac{e}{F}$. Thus,
\begin{eqnarray*}
  G_1 &=& (c_{111}g_1+c_{211}g_2+c_{311}g_3,c_{121}g_1+c_{221}g_2+c_{321}g_3,c_{131}g_1+c_{231}g_2+c_{331}g_3), \\
  G_2 &=& (c_{112}g_1+c_{212}g_2+c_{312}g_3,c_{122}g_1+c_{222}g_2+c_{322}g_3,c_{132}g_1+c_{232}g_2+c_{332}g_3), \\
  G_3 &=& (c_{113}g_1+c_{213}g_2+c_{313}g_3,c_{123}g_1+c_{223}g_2+c_{323}g_3,c_{133}g_1+c_{233}g_2+c_{333}g_3),
\end{eqnarray*}
Moreover, $\langle E_i,G_j\rangle=\delta_{ij}$, where
$\delta_{ij}$ denotes the Kronecker delta and $E_i=e_iF$ is the
$\mathbb A$ product of $e_i$ and $F$.
\eprop\noindent\textbf{Proof.} The first part is proven in
\cite{EL1} p. 198. Consider the function $H$ defined by
(\ref{isolo}), then
$$
e_i=H_*E_i=(\langle G_1,E_i\rangle,\langle G_2,E_i\rangle,\langle
G_3,E_i\rangle).
$$
Thus, the proof is finished. $\Box$

\subsection{Multivalued rectifications of algebrizable vector fields}

If $F$ is an $\mathbb A$-algebrizable vector field, the local
inverse of $H$ given by (\ref{isolo}) carries Euclidean geodesics
to geodesics of a certain Riemannian metric associated with $F$.
The function $H$ is well defined locally on $\Omega(F,\mathbb
A^*)$ and globally can be defined as a multivalued function on
$\Omega(F,\mathbb A^*)$, as we see in the following result.

\bcor\label{ecalpra} Let $F$ be an $\mathbb A$-algebrizable vector
field defined on an
open set $\Omega$.\\
1. For each point $w_0\in\Omega(F,\mathbb A^*)$, there exists a
simply connected neighborhood $N(w_0)$ contained in
$\Omega(F,\mathbb A^*)$ such that $H$ given by (\ref{isolo}) is
well defined on $N(w_0)$. This function depends only on $F$ and at an initial point $w_0$.\\
2. The function $H$ defined locally by (\ref{isolo}) can be
extended to a multivalued function $H:\Omega(F,\mathbb
A^*)\to\mathbb R^2$, which depends only on $F$ and at an initial
point $w_0$ in each connected component of $\Omega(F,\mathbb A^*)$.\\
3. If all the connected components of $\Omega(F,\mathbb A^*)$ are
simply connected, then $H$ defined locally by (\ref{isolo}) can be
extended to an univalued function $H:\Omega(F,\mathbb
A^*)\to\mathbb R^3$, which depends only on $F$ and at an initial
point $w_0$ in $\Omega(F,\mathbb A^*)$.

If we change the points $w_0$ in $\Omega(F,\mathbb A^*)$ used to
define $H$, then $H$ changes by a adding a constant on the
connected components. \ecor \noindent\textbf{Proof.} The proof is
the same as in the case of planar vector fields, which appears in
\cite{Dyg}. $\Box$

\bdefe Let $c\in\mathbb R^3$. We say that a multivalued function
$H:\Omega\to\mathbb R^3$ is \emph{multivalued rectification} of a
vector field $F$ if for every $w_0\in\Omega$ there exists a
neighborhood $N(w_0)$ such that for each of the branches
$k_\alpha=H_\alpha(w_0)$, $\alpha\in I$ (enumerates the branches),
there exists a neighborhood $V_\alpha$ and a local diffeomorphism
$H_\alpha:N(w_0)\to V_\alpha$ with $H_{\alpha*}F=c$. \edefe

In this work we denote the \emph{natural Riemannian metric}
$\lambda$ on $\mathbb R^3$ given by $\lambda_{ij}=\delta_{ij}$ for
$i,j\in\{1,2,3\}$, where $\delta_{ij}$ denotes the  Kronecker
delta.

The rectifying coordinates of vector fields are used by the
classical Lie methods in order to find their integral curves. In
this paper we use $H$ given in (\ref{isolo}) to define a
Riemannian metric $g$ on $\Omega(F,\mathbb A^*)$. Therefore, the
geodesics of $g$ are carried to Euclidean straight lines. The
following theorem gives multivalued rectifications for
algebrizable vector fields on terms of line integrals.

\bthm\label{hr2} Let $F$ be an $\mathbb A$-algebrizable
vector field defined on an open set $\Omega\subset\mathbb R^3$:\\
1. There exists a multivalued rectification $H$ of $\alpha F$ on
$\Omega(F,\mathbb A^*)$ given by
$$
H(w)=\int_\gamma\frac{e}{F}d\xi,
$$
under which $H_*(F)=e$, where $e$ is the
unit of $\mathbb A$.\\
2. There exists a Riemannian metric $g$ on $\Omega(\alpha
F,\mathbb A^*)$ defined by
$$
g=\frac{1}{\|e\|^2}H^*\lambda,
$$
where $\|e\|$ is the Euclidian norm of $e$, and $\lambda$ is the
natural Riemannian metric of $\mathbb R^3$. \ethm
\noindent\textbf{Proof.} The proof is the same as the case of
planar vector fields, which appears in \cite{Dyg}. $\Box$

\bcor\label{gdf} If $F$ is an $\mathbb A$-algebrizable vector
field defined on an open set $\Omega$, then the Rimennian metric
$g=\frac{1}{\|e\|^2}H^*\lambda$ makes $F$ geodesible on
$\Omega(F,\mathbb A^*)$. \ecor

\section{Algebrizable vector fields and their first
integrals}\label{camposfygs}

In this section, we give three propositions in which we associate
with each $\mathbb A$-algebrizable vector field $F$ three vector
fields $E_i=e_iF$ for $i=1,2,3$ as well as three gradient vector
fields $G_i$ for $i=1,2,3$. These vector fields are obtained by
using the products of the algebras $\mathbb
A_r^3(p_1,\cdots,p_6)$, $\mathbb A_{r,s}^3(p_1,p_2)$ and $\mathbb
A_{1,2,3}^3$; see Section \ref{s1}.

We show in this section that the potential functions of the vector
fields $G_i$ for $i=1,2,3$ determine two first integral for $F$
such that the integral curves of $F$ are locally contained in the
intersections of the level surfaces of these first integrals.

\subsection{Cases considered by Ward: algebras with unit $e_r$}\label{auer}

In this subsection the algebra $\mathbb A=\mathbb
A^3_r(p_1,\cdots,p_6)$ with unit $e_r$ defined by six real
parameters $p_1,\cdots,p_6$ is considered. The structure constants
for $\mathbb A$ are the following:
$$
\begin{array}{ccc}
  C_{r11}=1, & C_{r12}=0, & C_{r13}=0, \\
  C_{r21}=0, & C_{r22}=1, & C_{r23}=0,  \\
  C_{r31}=0, & C_{r32}=0, & C_{r33}=1,  \\
   C_{s11}=0, & C_{s12}=1, & C_{s13}=0, \\
  C_{s21}=p_7, & C_{s22}=p_1, & C_{s23}=p_2,  \\
  C_{s31}=p_8, & C_{s32}=p_3, & C_{s33}=p_4,  \\
    C_{t11}=0, & C_{t12}=0, & C_{t13}=1, \\
  C_{t21}=p_8, & C_{t22} =p_3, & C_{t23}=p_4,  \\
  C_{t31}=p_9, & C_{t32}=p_5, & C_{t33}=p_6,  \\
\end{array}%
$$
where $\{r,s,t\}=\{1,2,3\}$. The first fundamental representation
of $F=(f_1,f_2,f_3)$ is
$$
R(F)=\left(%
\begin{array}{ccc}
  f_1 & p_7f_2+p_8f_3 & p_8f_2+p_9f_3 \\
  f_2 & f_1+p_1f_2+p_3f_3 & p_3f_2+p_5f_3 \\
  f_3 & p_2f_2+p_4f_3 & f_1+p_4f_2+p_6f_3 \\
\end{array}%
\right).
$$

By using the equalities given in (\ref{ceq}) and the commutativity
equations (\ref{ecd}), the determinant of $R(F)$ is
\begin{eqnarray*}
   \det(R(F))&=& f_r^3+(p_1+p_4)f_r^2f_s+(p_3+p_6)f_r^2f_t+(p_1p_4-p_2p_3-p_7)f_rf_s^2\\
    & + & (p_1p_6-p_2p_5-2p_8)f_rf_sf_t +(p_3p_6-p_4p_5-p_9)f_rf_t^2\\& + &(p_2p_8-p_4p_7)f_s^3
    + (2p_2p_9-p_4p_8-p_6p_7)f_s^2f_t \\& + &(p_5p_7-p_6p_8-p_1p_9+ p_4p_9)f_sf_t^2
    + (p_5p_8-p_3p_9)f_t^3.
\end{eqnarray*}

\bprop\label{proprst} Let $\mathbb A$ be an algebra $\mathbb A
=\mathbb A^3_r(p_1,\cdots,p_6)$ defined by six real parameters
$p_1,\cdots,p_6$ and $\{r,s,t\}=\{1,2,3\}$. Consider an $\mathbb
A$-algebrizable vector field $F=(f_1,f_2,f_3)$ defined over an
open set $\Omega\subset\mathbb R^3$, and $E_i:=e_i F$ for
$i=1,2,3$. Thus,
\begin{eqnarray*}
  E_r &=& e_rF=F, \\
  E_s &=& e_sF=
  (p_7f_s+p_8f_t)e_r+(f_r+p_1f_s+p_3f_t)e_s+(p_2f_s+p_4f_t)e_t, \\
  E_t &=&
  e_tF=(p_8f_s+p_9f_t)e_r+(p_3f_s+p_5f_t)e_s+(f_r+p_4f_s+p_6f_t)e_t.
\end{eqnarray*}
In addition,
$$
\int_\gamma \frac{e}{F}d\xi=\left(\int_\gamma G_1\cdot
ds,\int_\gamma G_2\cdot ds,\int_\gamma G_3\cdot ds\right),
$$
where $\gamma$ is a differentiable path contained in
$\Omega(F,\mathbb A^*)$ and $G_r$, $G_s$, $G_t$ are the
conservative vector fields given by
\begin{eqnarray*}
  G_r &=& \frac{1}{\det(R(F))} [[f_r^2+(p_1+p_4)f_rf_s+(p_3+p_6)f_rf_t+(p_1p_4-p_2p_3)f_s^2\\
  &+&(p_1p_6-p_2p_5)f_sf_t+ (p_3p_6-p_4p_5)f_t^2]e_r\\
  &+& [-p_7f_rf_s-p_8f_rf_t+(p_2p_8-p_4p_7)f_s^2+
    (p_2p_9-p_6p_7)f_sf_t+(p_5p_7-p_3p_8)f_t^2]e_s\\ &+
    & [-p_8f_rf_s-p_9f_rf_t+(p_2p_9-p_4p_8)f_s^2+(p_5p_7-p_1p_9)f_sf_t+
    (p_5p_8-p_3p_9)f_t^2]e_t],
\end{eqnarray*}
\begin{eqnarray*}
  G_s &=& \frac{1}{\det(R(F))}[[-f_rf_s-p_4f_s^2-(p_6-p_3)f_sf_t+p_5f_t^2]e_r\\
  &+& [f_r^2+p_4f_rf_s+p_6f_rf_t- p_8f_sf_t-p_9f_t^2]e_s\\ &-& [p_3f_rf_s-p_8f_s^2+p_5f_rf_t-p_9f_sf_t]e_t],
\end{eqnarray*}
\begin{eqnarray*}
  G_t &=&  \frac{1}{\det(R(F))} [[p_2f_s^2-f_rf_t-(p_1-p_4)f_sf_t-p_3f_t^2]e_r\\
    &-& [p_2f_rf_s+p_4f_rf_t-p_7f_sf_t-p_8f_t^2]e_s\\ &+&
    [f_r^2+p_1f_rf_s+p_3f_rf_t-p_7f_s^2-
  p_8f_sf_t]e_t].
\end{eqnarray*}
Moreover, $\langle E_i,G_j\rangle=\delta_{ij}$. Thus, the
potential function $h_i$ of $G_i$ given by (\ref{prim10}) is a
first integral for $E_j$ for $i,j\in\{1,2,3\}$ and $i\neq j$.
\eprop \textbf{Proof.} We consider only the case when $r=1$,
$s=2$, and $t=3$. Let $\mathbb A=\mathbb A^3_1(p_1,\cdots,p_6)$.

From Proposition \ref{palvf} and the structure constants of
$\mathbb A$, the vector fields $G_i$ for $i=1,2,3$ are
\begin{eqnarray*}
  G_1 &=& (k_1,p_7k_2+p_8k_3,p_8k_2+p_9k_3), \\
  G_2 &=& (k_2,k_1+p_1k_2+p_3k_3,p_3k_2+p_5k_3), \\
  G_3 &=& (k_3,p_2k_2+p_4k_3,k_1+p_4k_2+p_6k_3),
\end{eqnarray*}
where $(k_1,k_2,k_3)=\frac{e}{F}$, and
\begin{eqnarray*}
  k_1 &=& \frac{1}{\det(R(F))} [f_1^2+(p_1+p_4)f_1f_2+(p_3+p_6)f_1f_3\\ &+ &(p_1p_4-p_2p_3)f_2^2
  + (p_1p_6-p_2p_5)f_2f_3+(p_3p_6-p_4p_5)f_3^2],\\
    k_2 &=& \frac{1}{\det(R(F))}[-f_1f_2-p_4f_2^2-(p_6-p_3)f_2f_3+p_5f_3^2], \\
    k_3 &=& \frac{1}{\det(R(F))}[-f_1f_3+p_2f_2^2-(p_1-p_4)f_2f_3-p_3f_3^2].
\end{eqnarray*}

By substituting the $k_i's$ for $i=1,2,3$ in the $G_j$, the
expression of $G_j$ given above can be obtained for $j=1,2,3$.

The last part of Proposition \ref{palvf} gives $\langle
E_i,G_j\rangle=\delta_{ij}$. $\Box$

In the following example, we consider the radial vector field
$F(w)=w$ seen as an $\mathbb A$-algebrizable vector field with
respect to the algebra $\mathbb A=\mathbb A^3_1(0,\cdots,0)$
defined in Section \ref{s1}.

\bejem\label{gpcr} Consider the three-dimensional vector field
$F(w)=(x_1,x_2,x_3)$ and $\mathbb A=\mathbb A^3_1(0,\cdots,0)$.
The vector fields $E_i's$ are given by $E_i=x_1e_i$ for $i=1,2,3$.

By Subsection \ref{auer}, we have that
$$
R(F)=\left(%
\begin{array}{ccc}
  x_1 & 0 & 0 \\
  x_2 & x_1 & 0 \\
  x_3 & 0 & x_1 \\
\end{array}%
\right),\qquad [R(F)]^{-1}=\left(%
\begin{array}{ccc}
  \frac{1}{x_1} & 0 & 0 \\
  \frac{-x_2}{x_1^2} & \frac{1}{x_1} & 0 \\
  \frac{-x_3}{x_1^2} & 0 & \frac{1}{x_1} \\
\end{array}%
\right),
$$
then
$[F(w)]^{-1}=(\frac{1}{x_1},\frac{-x_2}{x_1^2},\frac{-x_3}{x_1^2})$.
In addition,
\begin{eqnarray*}
  \int_\gamma \frac{e}{F}d\xi &=& \int
\left(\frac{x_1}{x_1^2},\frac{-y}{x_1^2},\frac{-z}{x_1^2}\right)(\gamma_1',\gamma_2',\gamma_3')dt\\
   &=& \int
\left(\frac{\gamma_1'}{x_1},-\frac{\gamma_1'x_2}{x_1^2}+\frac{\gamma_2'}{x_1},-\frac{\gamma_1'x_3}{x_1^2}+\frac{\gamma_3'}{x_1}\right)dt\\
   &=& \int
(\langle G_1,\gamma'\rangle,\langle G_2,\gamma'\rangle,\langle
G_3,\gamma'\rangle)dt,
\end{eqnarray*}
where the vector fields $G_i's$ for $i=1,2,3$ are given by
$$
G_1=\left(\frac{1}{x_1},0,0\right),\,\,\,\,\,\,\,G_2=\left(-\frac{x_2}{x_1^2},\frac{1}{x_1},
0\right),\,\,\,G_3=\left(-\frac{x_3}{x_1^2},0,\frac{1}{x_1}\right),
$$
and they are conservative.

Note that $\langle E_i, G_j\rangle=\delta_{ij}$, where
$\delta_{ij}$ is the Kronecker delta and $E_i=e_iF$ for
$i,j\in\{1,2,3\}$.

In this case $\Omega(F,\mathbb A^*)=\{(x_1,x_2,x_3)\in\mathbb
R^3\,:\,x_1\neq 0\}$. \eejem

\subsection{Cases do not considered by Ward}\label{cncbw}

\subsubsection{Algebras with unit $e_r+e_s$}

Now the algebra $\mathbb A=\mathbb A^3_{r,s}(p_1,p_2)$ defined by
two real parameters $p_1$ and $p_2$ is considered. The structure
constants are given by
$$
\begin{array}{ccc}
  C_{rrr}=1, & C_{rrs}=0, & C_{rrt}=0, \\
  C_{rsr}=0, & C_{rss}=0, & C_{rst}=0,  \\
  C_{rtr}=0, & C_{rts}=0, & C_{rtt}=0,  \\
   C_{srr}=0, & C_{srs}=0, & C_{srt}=0, \\
  C_{ssr}=0, & C_{sss}=1, & C_{sst}=0,  \\
  C_{str}=0, & C_{sts}=0, & C_{stt}=1,  \\
    C_{trr}=0, & C_{trs}=0, & C_{trt}=0, \\
  C_{tsr}=0, & C_{tss} =0, & C_{tst}=1,  \\
  C_{ttr}=0, & C_{tts}=p_1, & C_{ttt}=p_2.  \\
\end{array}%
$$

The determinant of the first fundamental representation of
$F=f_re_r+f_se_s+f_te_t$ is
$\det(R(F))=f_rf_s^2+p_2f_rf_sf_t-p_1f_rf_t^2$.

\bprop\label{ancbw} Let $\mathbb A$ be an algebra defined by the
two real parameters $p_1 $, $p_2$ and $\{r,s,t\}=\{1,2,3\}$.
Consider an $\mathbb A$-algebrizable vector field
$F=(f_1,f_2,f_3)$ defined on an open set $\Omega\subset\mathbb
R^3$, and $E_i:=e_i F$ for $i=1,2,3$. Thus,
\begin{eqnarray*}
  E_r &=& e_rF=f_re_r, \\
  E_s &=& e_sF=f_se_s+f_te_t, \\
  E_t &=& e_tF=p_1f_te_s+(f_s+p_2f_t)e_t.
\end{eqnarray*}
In addition,
$$
\int_\gamma \frac{e}{F}d\xi=\left(\int_\gamma G_1\cdot
ds,\int_\gamma G_2\cdot ds,\int_\gamma G_3\cdot ds\right),
$$
where $\gamma$ is a differentiable path contained in
$\Omega(F,\mathbb A^*)$ and $G_r$, $G_s$, $G_t$ are the
conservative vector fields given by
\begin{eqnarray*}
  G_r &=&
  \frac{1}{f_r}e_r,\\
  G_s &=& \frac{f_s+p_2f_t}{f_s^2+p_2f_sf_t-p_1f_t^2}e_s-\frac{p_1f_t}{f_s^2+p_2f_sf_t-p_1f_t^2}e_t,\\
  G_t &=&
  \frac{-f_t}{f_s^2+p_2f_sf_t-p_1f_t^2}e_s+\frac{f_s}{f_s^2+p_2f_sf_t-p_1f_t^2}e_t.
\end{eqnarray*}
Moreover, $\langle E_i,G_j\rangle=\delta_{ij}$. Thus, the
potential function $h_i$ of $G_i$ given by (\ref{prim10}) is a
first integral for $E_j$ for $i,j\in\{1,2,3\}$ and $i\neq j$.
\eprop \textbf{Proof.} The multiplicative inverse
$\frac{e}{F}=(k_1,k_2,k_3)$ of $F$ is given by
$$
\frac{e}{F}=\left(\frac{1}{f_r}e_r+\frac{f_s+p_2f_t}{f_s^2+p_2f_sf_t-p_1f_t^2}e_s-\frac{f_t}{f_s^2+p_2f_sf_t-p_1f_t^2}e_t\right),
$$
then
$$
k_r=\frac{1}{f_r},\qquad
k_s=\frac{f_s+p_2f_t}{f_s^2+p_2f_sf_t-p_1f_t^2},\qquad
k_t=\frac{-f_t}{f_s^2+p_2f_sf_t-p_1f_t^2}.
$$

From Proposition \ref{palvf} and the structure constants of
$\mathbb A$, the vector fields $G_i$ for $i=1,2,3$ are
$$
G_r = k_re_r,\quad  G_s = k_se_s+p_1k_te_t,\quad G_t =
k_te_s+(k_s+p_2k_t)e_t.
$$

By substituting the $k_i's$ for $i=1,2,3$ in $G_j$, the expression
of $G_j$ given above can be obtained for $j=1,2,3$.

The last part of Proposition \ref{palvf} gives $\langle
E_i,G_j\rangle=\delta_{ij}$.$\Box$

\subsubsection{Algebra with unit $e=e_1+e_2+e_3$}

Now the algebra $\mathbb A=\mathbb A^3_{1,2,3}$ is considered. The
structure constants $c_{i,j,k}$ for $i,j,k\in\{1,2,3\}$ are given
by
$$
c_{ijk}=
\left\{%
\begin{array}{ll}
   1 , & \hbox{$i=j=k$;} \\
   0, & \hbox{other case.} \\
\end{array}%
\right.
$$

The determinant of the first fundamental representation of
$F=f_1e_1+f_2e_2+f_3e_3$ is $\det(R(F))=f_1f_2f_3$.

\bprop\label{aecu123} Let $\mathbb A$ be the algebra $\mathbb
A^3_{1,2,3}$. Consider an $\mathbb A$-algebrizable vector field
$F=(f_1,f_2,f_3)$ defined over an open set $\Omega\subset\mathbb
R^3$, and $E_i:=e_i F$ for $i=1,2,3$. Thus,
$$
F_1=(f_1,0,0),\quad F_2=(0,f_2,0), \quad F_2=(0,0,f_3).
$$
In addition,
$$
\int_\gamma \frac{e}{F}d\xi=\left(\int_\gamma G_1\cdot
ds,\int_\gamma G_2\cdot ds,\int_\gamma G_3\cdot ds\right),
$$
where $\gamma$ is a differentiable path contained in
$\Omega(F,\mathbb A^*)$ and $G_1$, $G_2$, $G_3$ are the
conservative vector fields given by $G_i=\frac{1}{f_i}e_i$ for
$i=1,2,3$. Moreover, $\langle E_i,G_j\rangle=\delta_{ij}$. Thus,
the potential function $h_i$ of $G_i$ given by (\ref{prim10}) is a
first integral for $E_j$ for $i,j\in\{1,2,3\}$ and $i\neq j$.
\eprop \textbf{Proof.} The proof is trivial. $\Box$

\subsection{The regular domain on $F$}\label{rdof}

By using the first fundamental representation $R$ of $\mathbb A$,
we notice that $F(w)$ is regular in $\mathbb A$ if and only if
$\det(R(F(w)))\neq 0$. Thus, from Subsections \ref{auer} and
\ref{cncbw} we obtain that $F(w)$ is regular for $w$ in the
corresponding set
\begin{enumerate}
    \item $\Omega(F,\mathbb A^*)=P^{-1}(\mathbb R^*)$ if $F$ is $\mathbb
    A$-algebrizable for $\mathbb
    A=A^3_{r}(p_1,\cdots,p_6)$, where $P$ is the polynomial function
\begin{eqnarray*}
  P &=& f_r^3+(p_1+p_4)f_r^2f_s+(p_3+p_6)f_r^2f_t+(p_1p_4-p_2p_3-p_7)f_rf_s^2\\
    & + & (p_1p_6-p_2p_5-2p_8)f_rf_sf_t +(p_3p_6-p_4p_5-p_9)f_rf_t^2\\& + &(p_2p_8-p_4p_7)f_s^3
    + (2p_2p_9-p_4p_8-p_6p_7)f_s^2f_t \\& + &(p_5p_7-p_6p_8-p_1p_9+ p_4p_9)f_sf_t^2
    + (p_5p_8-p_3p_9)f_t^3,
\end{eqnarray*}
and $\{r,s,t\}=\{1,2,3\}$.
    \item $\Omega(F,\mathbb A^*)=(f_rf_s^2+p_2f_rf_sf_t-p_1f_rf_t^2)^{-1}(\mathbb R^*)$ if $F$ is $\mathbb
    A$-algebrizable for $\mathbb
    A=A^3_{r,s}(p_1,p_2)$, where $\{r,s,t\}=\{1,2,3\}$.
    \item $\Omega(F,\mathbb A^*)=(f_1f_2f_3)^{-1}(\mathbb R^*)$ if $F$ is $\mathbb
    A$-algebrizable for $\mathbb
    A=A^3_{1,2,3}$.
\end{enumerate}
Thus, generically the set $\Omega(F,\mathbb A^*)$ is dense in
$\Omega$.

\subsection{First integrals of algebrizable three-dimensional vector fields}\label{camposfyg}

The level surfaces of a first integral contain the integral curves
of the corresponding vector field. If $G=(p,q,r)$ is a gradient
vector field and $F$ is a vector field with inner product $\langle
F,G\rangle=0$, then a first integral $h$ of $F$ may be calculated
by quadrature in the form
\begin{equation}\label{prim10}
h(w)=\int p(w) dx_1+\int\left(q(w)-\frac{\partial}{\partial
x_2}\int p(w)dx_1\right)dx_2,
\end{equation}
where $w=(x_1,x_2,x_3)$.

\bdefe We say that a function $H$ is an $\mathbb
A$-\textit{antiderivative} of an $\mathbb A$-algebrizable vector
field $F$ if the $\mathbb A$-derivative of $H$ is $H'=F$. \edefe

\bdefe We say that two differentiable functions
$h_1,h_2:\Omega\subset\mathbb R^3\to\mathbb R$ defined on an open
set $\Omega\subset\mathbb R^3$ have \emph{transversal
intersections}, if for every pair of level surfaces
$S_1=h_1^{-1}(c_1)$, $S_2=h_2^{-1}(c_2)$ with $S_1\cap
S_2\neq\emptyset$, $S_1$ and $S_2$ have transversal intersection.
\edefe

The following proposition gives two first integrals having
transversal intersections for algebrizable three-dimensional
vector fields.

\bprop\label{calprin} Let $F$ be an $\mathbb A$-algebrizable
vector field defined on an open set $\Omega\subset\mathbb R^3$,
$w_0\in\Omega(F,\mathbb A^*)$, and $N(w_0)$ a simply connected
open set contained in $\Omega(F,\mathbb A^*)$. A function
$H=(h_1,h_2,h_3)$ is an $\mathbb A$-antiderivative of
$\frac{e}{F}$ on $N(w_0)$ if and only if for $w_0,w\in N(w_0)$
$$
H(w)=H(w_0)+\left(\int_{w_0}^w G_1\cdot ds,\int_{w_0}^w  G_2\cdot
ds,\int_{w_0}^w  G_3\cdot ds\right).
$$
If $H$ is an $\mathbb A$-antiderivative of $\frac{e}{F}$ on
$N(w_0)$, then a pair of first integrals for $F$ on $N(w_0)$ is
given as follows:
\begin{enumerate}
\item $h_s$ and $h_t$, when $F$ is $\mathbb A$-algebrizable with respect to $\mathbb A=\mathbb
A^3_r(p_1,\cdots,p_6)$;
\item $h_r-h_s$ and $h_r-h_s+h_t$, when $F$ is $\mathbb A$-algebrizable with respect to $\mathbb A=\mathbb A^3_{r,s}(p_1,p_2)$; and
\item $h_r-h_t$ and $h_r-h_s$, when $F$ is $\mathbb A$-algebrizable with respect to $\mathbb A=\mathbb
A^3_{1,2,3}$,
\end{enumerate}
where $\{r,s,t\}=\{1,2,3\}$. Moreover, these pairs of first
integrals have transversal intersections. \eprop
\noindent\textbf{Proof.} By Proposition \ref{palvf}, there exist
conservative vector fields $G_1$, $G_2$, and $G_3$ satisfying
(\ref{g12}). Let
$$
K(w):=\int_{w_0}^w \frac{e}{F}dw.
$$
If $H=(h_1,h_2,h_3)$ is an $\mathbb A$-antiderivative of
$\frac{e}{F}$ on $N(w_0)$, then $K(w)=H(w)+C$, where $C\in\mathbb
A$ is a constant. Thus, $K(w)=H(w)-H(w_0)$. Now if $H$ has the
expression given in the proposition, then by the fundamental
theorem of calculus on algebras, $H$ is an antiderivative of
$\frac{e}{F}$.

From Proposition \ref{proprst}, we have $e=e_r$, then $F=E_r$. By
Propositions \ref{proprst}, \ref{ancbw}, and \ref{aecu123}, the
vector fields $E_i,G_i$ for $i=1,2,3$ satisfy $\langle
E_i,G_j\rangle=\delta_{ij}$ where $\delta_{ij}$ denotes the
Kronecker delta. Thus, the potential functions $h_s$, $h_t$ of the
gradient vector fields $G_s$, $G_t$ defined on $\Omega(F,\mathbb
A^*)$ are first integrals of $E_r$. Therefore, item 1) is proved.

If $H=(h_1,h_2,h_3)$ is an $\mathbb A$-antiderivative of
$\frac{e}{F}$, then $\langle\nabla h ,F\rangle=0$ for each
candidate $h$ to be a first integral of $F$. Thus, items 2) and 3)
can be proved.

Suppose the existence of a point $w_0\in\Omega(F,\mathbb A^*)$
such that the intersection of the level surfaces of $h_s$, $h_t$
by $w_0$ does not have a transversal intersection, then the
tangent spaces at $w_0$ coincide and the vector fields $G_s$,
$G_t$ are linearly dependent. That is, there exists $c\in\mathbb
R$ such that $G_s+cG_t=0$, and by Proposition \ref{proprst}
$\langle E_i,G_j\rangle=\delta_{ij}$ for $i,j\in\{1,2,3\}$, where
$\delta_{ij}$ denotes the Kronecker delta. Thus,
$$
0=\langle E_s,G_s+cG_t\rangle=\langle E_s,G_s\rangle+c\langle
E_s,G_t\rangle=1,
$$
which is a contradiction. As a result $h_s$ and $h_t$ have a
transversal intersection at $w_0$. In the same way, since the
vectorial products following $(G_r-G_s)\times(G_r-G_s+G_t)\neq 0$
and $(G_r-G_t)\times (G_r-G_s)\neq 0$, the transversal
intersections for the the cases of items 2) and 3) can be
obtained. $\Box$

The intersections of the level surfaces of two first integrals
having transversal intersections of a three-dimensional vector
field $F$ define solutions of $F$, in the sense of the following
corollary.

\bcor\label{titfi} Let $h_1,h_2:\Omega\in\mathbb R^3\to\mathbb R$
be first integrals having transversal intersections of a vector
field $F$ defined on an open set $\Omega\in\mathbb R^3$. If
$w_0\in\Omega$ and $\mathcal L^1_{w_0}$, $\mathcal L^2_{w_0}$ are
level surfaces of $h_1$, $h_2$ containing $w_0$. Thus, for each
$w_1\in \mathcal L^1_{w_0}\cap \mathcal L^2_{w_0}$, there exists
an open interval $I\subset\mathbb R$ containing $0$ such that
$$
\{\Phi_t(w_1)\,\,:\,\,t\in I\}\subset \mathcal L^1_{w_0}\cap
\mathcal L^2_{w_0},
$$
where $\Phi_t(\cdot)$ is the flow of $F$. \ecor

The $\mathbb A$-derivative satisfies the usual differentiation
rules of the complex derivative. For example, the antiderivatives
of the polynomial functions in the variable of $\mathbb A$, such
as
$$w\mapsto a_0+a_1w+a_2w^2+\cdots+a_nw^n$$
where $a_0,a_1,\cdots,a_n\in\mathbb A$ are constants, can be
expressed by
$$H(w)=c+a_0w+\frac{1}{2}a_1w^2+\frac{1}{3}a_2w^3+\cdots+\frac{1}{n+1}a_nw^{n+1}.$$
Thus, for an important family of vector fields $F$, the first
integrals can be easily calculated by their antiderivatives by
using Corollary \ref{calprin}. The $\mathbb A$-algebrizable vector
fields of the form $F(w)=w^n$ (for $n\neq 1$) have first integrals
which can be calculated by finding antiderivatives of the
functions $e/F(w)=w^{-n}$, as we see in the following example.

\bejem\label{pipp} Let $F$ be an $\mathbb A$-algebrizable vector
field $F(w)=w^n$ for $n\neq 1$, defined on
$$\Omega=
\left\{%
\begin{array}{ll}
    \mathbb A^*, & \hbox{if $n<1$;} \\
    \mathbb R^2, & \hbox{if $n>0$;} \\
\end{array}%
\right.
$$
therefore, $H(w)=\frac{-e}{(n-1)w^{n-1}}$ is an antiderivative of
$\frac{e}{F}$ defined on $\Omega(F,\mathbb A^*)$. Thus, if
$$(h_1(w),h_2(w),h_3(w))=\frac{-e}{(n-1)w^{n-1}},$$ then $h_i$ is a first
integral for $E_j=e_jF$ and $E_k=e_kF$ on $\Omega(F,\mathbb A^*)$
for $\{i,j,k\}=\{1,2,3\}$. \eejem

\section{Riemannian metrics and commutative orthonormal
frames}\label{rmcof}

\subsection{On Riemannian metrics and orthonormal frames}\label{mgc}

The set of vector fields $\{E_1,E_2,E_3\}$ defines an orthonormal
frame, as we see in the following proposition.

\bprop\label{ortfr} Let $F$ be a vector field defined on an open
set $\Omega\subset\mathbb R^3$ that is $\mathbb A$-algebrizable
with respect to an algebra $\mathbb A$, and $E_i$ the vector field
$E_i=e_iF$ for $i=1,2,3$, where $e_iF$ denotes the product of
$e_i$ and $F$ with respect to $\mathbb A$. Thus, the set of vector
fields $\{E_1,E_2,E_3\}$ defines an orthonormal frame on
$\Omega(F,\mathbb A^*)$ for the Riemannian metric
$g=\frac{1}{\|e\|^2}H^*\lambda$ defined in Theorem \ref{hr2}.
\eprop \textbf{Proof.} We have that $H_*E_i=e_i$ for $i=1,2,3$.
Since $g=\frac{1}{\|e\|^2}H^*\lambda$, we have that
$\{E_1,E_2,E_3\}$ is orthonormal with respect to $g$.$\Box$

\bprop\label{tenmr} Let $F$ be a vector field defined on an open
set $\Omega\subset\mathbb R^3$ that is $\mathbb A$-algebrizable
with respect to an algebra $\mathbb A$, and $E_i$ the vector field
$E_i=e_iF$ for $i=1,2,3$, where $e_iF$ denotes the product of
$e_i$ and $F$ with respect to $\mathbb A$. Thus, the metric
$g=\frac{1}{\|e\|^2}H^*\lambda$ defined on $\Omega(F,\mathbb A^*)$
is given by
$$
g=\left(%
\begin{array}{ccc}
  a_1 & a_2 & a_3 \\
  a_2 & a_4 & a_5 \\
  a_3 & a_5 & a_6 \\
\end{array}%
\right),
$$
where $a_1$, $a_2$, $a_3$, $a_4$, $a_5$, and $a_6$ are the
solutions of the linear system of the six equations:
\begin{equation}\label{eclin}
    E_i g E_i^T=1,\quad i=1,2,3, \quad E_i g E_j^T=0, \quad
i,j\in\{1,2,3\},\quad i< j.
\end{equation}
\eprop \textbf{Proof.} By Proposition \ref{ortfr}, we have that
equations (\ref{eclin}) are satisfied. $\Box$

\bcor\label{c1mr} If $F=(f_1,f_2,f_2)$ is $\mathbb A$-algebrizable
with respect to the algebra $\mathbb A=\mathbb A^3_1(0,\cdots,0)$,
then $g$ defined on $\Omega(F,\mathbb A^*)=\{w\in\mathbb
R^3\,:\,f_1(w)\neq 0\}$ is given by
$$
g=\frac{1}{f_1^4}\left(%
\begin{array}{ccc}
  f_1^4+2f_2^2+f_3^2 & -f_1f_2 & -f_1f_3 \\
  -f_1f_2 & f_1^2 & 0 \\
  -f_1f_3 & 0 & f_1^2 \\
\end{array}%
\right).
$$
\ecor\textbf{Proof.} The vector fields $E_i$ are given in
Proposition \ref{proprst}. $\Box$

\bcor\label{c2mr} If $F=(f_1,f_2,f_2)$ is $\mathbb A$-algebrizable
with respect to $\mathbb A=\mathbb A^3_{1,2}(p_1,p_2)$, then $g$
defined on
$$
\Omega(F,\mathbb A^*)=\{w\in\mathbb R^3\,:\,f_1(w)([f_2(w)]^2 -p_1
f_2(w) f_3(w) - p_2 [f_3(w)]^2)\neq 0\}
$$
is given by
$$
g=\left(%
\begin{array}{ccc}
  \frac{1}{f_1^2} & 0 & 0 \\
  0 & \frac{(f_2-p_1f_3)^2 + f_3^2 }{(f_2^2 - p_1 f_2 f_3 - p_2 f_3^2)^2} & \frac{(-f_2 - p_2 f_2 + p_1p_2 f_3)f_3 }{(f_2^2 - p_1 f_2 f_3 - p_2 f_3^2)^2} \\
  0 & \frac{(-f_2 - p_2 f_2 + p_1p_2 f_3)f_3 }{(f_2^2 - p_1 f_2 f_3 - p_2 f_3^2)^2}
  & \frac{f_2^2 + p_2^2 f_3^2}{(f_2^2 - p_1 f_2 f_3 - p_2 f_3^2)^2} \\
\end{array}%
\right).
$$
\ecor\textbf{Proof.} The vector fields $E_i$ are given in
Proposition \ref{ancbw}. We use the results given in \cite{Dyg}.
$\Box$

\bcor\label{c3mr} If $F=(f_1,f_2,f_2)$ is $\mathbb A$-algebrizable
with respect to the algebra $\mathbb A=\mathbb A^3_{1,2,3}$, the
Riemannian metric $g$ defined on
$$\Omega(F,\mathbb A^*)=\{w\in\mathbb
R^3\,:\,f_1(w)f_2(w)f_3(w)\neq 0\}$$ by is given by
$g_{ij}=\frac{\delta_{ij}}{f_i^2}$ for $i,j\in\{1,2,3\}$, where
$\delta_{ij}$ is the Kronecker delta. \ecor\textbf{Proof.} The
vector fields $E_i$ are given in Proposition \ref{aecu123}. $\Box$

\subsection{Families of geodesible vector fields}

Let $\mathbb A$ be a three-dimensional algebra, $F=(f_1,f_2,f_3)$
an $\mathbb A$-algebrizable vector field, and $F_{\theta,\phi}$
the vector field obtained from the product
$F_{\theta,\phi}=u_{\theta,\phi}F$ with respect to $\mathbb A$ of
the unitary vector
$$
u_{\theta,\phi}=(\cos\theta\sin\phi,\sin\theta\sin\phi,\cos\phi)
$$
and $F$. Thus, we define the family of three-dimensional vector
fields
$$
\mathcal F=\{F_{\theta,\phi}\,:\,0\leq\theta\leq
2\pi,\,0\leq\phi\leq \pi\}.
$$

The vector fields of $\mathcal F$ are geodesible and commutative
with respect to the Lie bracket, as we see in the following
proposition.

\bprop\label{fcavf} Let $\mathbb A$ be a three-dimensional algebra
and $F=(f_1,f_2,f_3)$ an $\mathbb A$-algebrizable vector field.
The vector fields in $\mathcal F$ are geodesible with respect to
the Riemannian metric $g=\frac{1}{\|e\|^2}H^*\lambda$ given in
Proposition \ref{hr2}, and
$\{E_1,E_2,E_3\}=\{F_{0,\pi/2},F_{\pi/2,\pi/2},F_{0,0}\}$.
Moreover, $\mathcal F$ is commutative with respect to the Lie
bracket. \eprop \textbf{Proof.} We have that $H_*F_{\theta,\phi}$
is the unitary vector $u_{\theta,\phi}$, then by the definition of
$g$ the integral curves of $F_{\theta,\phi}$ are unitary geodesics
of $g$.

For $k,l\in\mathbb A$ we take $K_1=kF$ and $K_2=lF$. The Lie
bracket $[K_1,K_2]$ of every pair of vector fields $K_1$ and $K_1$
satisfies
$$[K_1,K_2]=dK_2(w)(K_1(w))-dK_1(w)(K_2(w)).$$
Since $K_1$ and $K_2$ are $\mathbb A$-algebrizable vector fields,
we have
\begin{eqnarray*}
  dK_2(w)(K_1(w))-dK_1(w)(K_2(w)) &=& K'_2(w)K_1(w)-K'_1(w)K_2(w) \\
   &=& l(F)'(w)kF(w)-k(F)'(w)lF(w)\\
   &=& 0,
\end{eqnarray*}
where $d$ denotes the usual differentiation and $'$ denotes the
$\mathbb A$-differentiation. Therefore, the vector fields $K_1$
and $K_2$ commute.$\Box$

\subsection{Metrics defined by Riemannian metrics}\label{med}

Let us remember that $\|\cdot\|$ denotes the Euclidean norm in
$\mathbb R^3$.

The following theorem relates the metric $d_F$ arising from the
Riemannian metric $g$ given in Theorem \ref{hr2} with the
antiderivatives $H$ of $\frac{e}{F}$ given in (\ref{isolo}).

\bthm\label{mh} Let $F$ be an $\mathbb A$-algebrizable vector
field defined on an open set $\Omega\subset\mathbb R^3$,
$N\subset\Omega(F,\mathbb A^*)$ a simply connected neighborhood,
and $H$ an antiderivative of $\frac{e}{F}$. If $N$ is small
enough, then
$$
d_{F}(w,w')=\|H(w)-H(w')\|
$$
is the metric arising from the Riemannian metric $g$ given in
Theorem \ref{hr2}. \ethm \noindent\textbf{Proof.} By Proposition
\ref{ortfr}, $\{E_1,E_2,E_3\}$ defines an orthonormal frame for
$g$ on $\Omega(F,\mathbb A^*)$. Proposition \ref{palvf} applied to
$\frac{e}{F}$ implies for trajectories $\beta_i(t)$ of $E_i=e_iF$
that
$$
h_i(w)-h_i(w')=\int_{\beta_i} G_i\cdot ds=\int_0^{t_i}\langle G_i,
E_i\rangle dt=\int_0^{t_i}dt
$$
is the $g$ length of  $\beta_i$ for $i=1,2,3$, where
$H=(h_1,h_2,h_3)$ is an antiderivative of $\frac{e}{F}$ on $N$;
see Proposition \ref{calprin}. Thus,
$$
\|H(w)-H(w')\|=\sqrt{\sum_{i=1}^3(h_i(w)-h_i(w'))^2}.\qquad\Box
$$

In order to show how the algebra structure can be used to obtain
the metric $d_F$ of Theorem \ref{mh} associated with the
Riemannian metric $g$ given in Theorem \ref{hr2}, we give the
following example which is closely related to Example \ref{pipp}.

\bejem Let $N\subset\mathbb A^*$ be a simply connected
neighborhood and $F$ a three-dimensional vector field having the
form $F(w)=w^n$ for $n\neq 1$ with respect to $\mathbb A$. If $N$
is small enough, then the metric $d_{F}$ defined by the Riemannian
metric $g$ given in Theorem \ref{hr2} satisfies
$$
d_F(w,w_0)=\left\|\frac{e}{(n-1)w^{n-1}}-\frac{e}{(n-1)w_0^{n-1}}\right\|.
$$
\eejem

Theorem \ref{mh} has the following interpretation. Let us define
$d_{F}:N\times N\to\mathbb R$ by
$$d_{F}(w,w')=\left\|\left(\int_w^{w'}G_1\cdot ds,\int_w^{w'} G_2\cdot
ds,\int_w^{w'} G_3\cdot ds\right)\right\|,$$ where $G_1$, $G_2$,
and $G_3$ are the conservative vector fields given in Proposition
\ref{palvf} applied to $\frac{e}{F}$. The vector fields $E_i$ are
defined by the product $E_i=e_iF$ for $i=1,2,3$. By Propositions
\ref{ortfr} and \ref{fcavf} the vectors $E_1$, $E_2$, and $E_3$
are perpendicular with respect to $g$ and commute with respect to
the Lie bracket. Thus, if $\{E_1,E_2,E_3\}=\{E_i,E_j,E_k\}$, then
for each point $w\in N$ the vector fields $\{E_j,E_k\}$ locally
define a family of surfaces $\mathcal L^i_w$, which is invariant
under the flow $\Phi^i_{t}$ of $E_i$. If we consider $w,w'\in N$,
then there exist two points $w_1,w_2$ and times $t_i$ for
$i=1,2,3$ such that $\Phi^i_{t_1}(\mathcal L^i_w)=\mathcal
L^i_{w_1}$, $\Phi^j_{t_2}(\mathcal L^j_{w_1})=\mathcal L^j_{w_2}$,
and $\Phi^k_{t_3}(\mathcal L^k_{w_2})=\mathcal L^k_{w_3}$. Under
these conditions, we have that
$d_{F}(w,w')=\sqrt{t_1^2+t_2^2+t_3^2}$.

In the case of planar vector fields, the above interpretation is
well-known for $\mathbb A$-algebrizable vector fields; see
\cite{Dyg}. For the case of the complex numbers, see \cite{Mur}.

\bejem Consider the algebra $\mathbb A=\mathbb A^3_1(0,\cdots,0)$
and the vector field $F=(x_1^2,2x_1x_2,2x_2x_3)$. Thus, $F$ is
$\mathbb A$-algebrizable and $F(w)=w^2$ with respect to the
variable of $\mathbb A$. By Corollary \ref{c1mr} the Riemannian
metric $g$ is defined on $\Omega(F,\mathbb
A^*)=\{(x_1,x_2,x_3)\in\mathbb R^3\,:\,x_1\neq 0\}$, and it is
given by
$$
g=\frac{1}{x_1^8}\left(%
\begin{array}{ccc}
  x_1^8+8x_1^2x_2^2+8x_1^3x_2^3 & -2x_1^3x_2 & -2x_1^3x_3 \\
  -2x_1^3x_2 & x_1^4 & 0 \\
  -2x_1^3x_3  & 0 & x_1^4\\
\end{array}%
\right).
$$
The set of regular points $\Omega(F,\mathbb A^*)$ of $F$ is given
by $\Omega(F,\mathbb A^*)=L^+\cup L^-$, where
$$L^+=\{(x_1,x_2,x_3)\in\mathbb
R^3\,:\,x_1>0\},\qquad L^-=\{(x_1,x_2,x_3)\in\mathbb
R^3\,:\,x_1<0\}.$$

By Example \ref{pipp}, $H(w)=\frac{-e}{w}$ is an antiderivative of
$\frac{e}{F}$, then by Example \ref{gpcr}
$H(x_1,x_2,x_3)=\left(-\frac{1}{x_1},\frac{x_2}{x_1^2},\frac{x_3}{x_1^2}\right)$.
We have that $H(w)=(h_1(w),h_2(w),h_3(w))$. Thus, if
$w=(x_1,x_2,x_3)$ and $w'=(x_1',x_2',x_3')$ are $L^+$ (or in
$L^-$), then
$$
d_F(w,w')=\sqrt{\left(-\frac{1}{x_1}+\frac{1}{x_1'}\right)^2
+\left(\frac{x_2}{x_1^2}-\frac{x_2'}{x_1'^2}\right)^2+\left(\frac{x_3}{x_1^2}-\frac{x_3'}{x_1'^2}\right)^2}.
$$

Moreover, by Corollary \ref{titfi} the level surfaces $\mathcal
L^s_{w_0}$, $\mathcal L^t_{w_0}$ of $h_s$, $h_t$ by $w_0$ contain
the integral curve of $E_r$ by $w_0$ and have transversal
intersections, then this integral curve locally coincides with the
intersection $S_s\cap S_t$.

If $w_0=(a,b,c)$ with $a\neq 0$, then $\mathcal L^2_{w_0}$ is
given by $\frac{x_2}{x_1^2}=\frac{b}{a^2}$ and $\mathcal
L^3_{w_0}$ by $\frac{x_3}{x_1^2}=\frac{c}{a^2}$. Thus, the curve
$\{\left(t,\frac{b}{a^2}t^2,\frac{c}{a^2}t^2\right)\,:\;t\in I\}$
for $I=(0,\infty)$ or $I=(-\infty,0)$ is a geodesic of $g$. \eejem

\end{document}